\documentclass[12pt]{article}

\usepackage{amsmath}
\usepackage{amssymb}
\usepackage{amsthm}
\usepackage{a4wide}
\usepackage[english]{babel}
\usepackage{stmaryrd} 

\numberwithin{equation}{section}

\newtheorem{Th}{Theorem}[section]
\newtheorem{Lemma}[Th]{Lemma}
\newtheorem{Prop}[Th]{Proposition}
\newtheorem{Cor}[Th]{Corollary}

\newtheorem{Rem}[Th]{Remark}

\newcommand{\Dem}{\noindent{\bf Proof }}

\newcommand{\N}{\mathbb{N}}
\newcommand{\Z}{\mathbb{Z}}
\newcommand{\Q}{\mathbb{Q}}
\newcommand{\R}{\mathbb{R}}
\newcommand{\C}{\mathbb{C}}
\newcommand{\K}{\mathbb{K}}

\newcommand{\eps}{\varepsilon}

\newcommand{\Span}{{\rm Span}}

\newcommand{\unk}{\{1,\ldots,k\}}

\newcommand{\unp}{\{1,\ldots,p\}}

\newcommand{\unq}{\{1,\ldots,q\}}

\newcommand{\cale}{{\mathcal{E}}}

\newcommand{\combitiny}[2]{{\tiny \left( \! \!  \begin{array}{c} #1 \\ #2 \end{array} \! \! \right)}}
\newcommand{\combi}[2]{{  \left( \begin{array}{c} #1 \\ #2 \end{array} \right)}}

\newcommand{\pha}{\varphi}
\renewcommand{\Re}{{\rm Re}}
\renewcommand{\Im}{{\rm Im}}

\newcommand{\rk}{{\rm rk}}

\newcommand{\om}{\omega}

\newcommand{\alrab}{\alpha_{r,a,b}}
\newcommand{\Qrab}{Q_{r,a,b}}
\newcommand{\Iab}{I_{ a,b}}
\newcommand{\chirab}{\chi_{n}}

\newcommand{\epsral}{\eps_{ \lambda}}

\newcommand{\omral}{\om_{  \lambda}}
\newcommand{\phiral}{\pha_{  \lambda}}

\newcommand{\matc}{[\clb]_{ \lambda, \beta\in\cale}}
\newcommand{\matd}{[\dbl]_{\beta,\lambda\in\cale}}
\newcommand{\dbl}{d_{\beta,\lambda}^{(b)}}
\newcommand{\clb}{c_{\lambda, \beta}^{(b)}}

\newcommand{\jln}{J_{\lambda,n}}

\newcommand{\sln}{S_{\lambda,n}}
\newcommand{\ibn}{I_{\beta,n}}

\newcommand{\itibn}{\widetilde I_{\beta,n}}
\newcommand{\ltibn}{\widetilde{\ell}_{\beta,n}}
\newcommand{\lti}{\widetilde{\ell}}

\newcommand{\ro}{\varrho}

\newcommand{\Rab}{{\mathcal R}_{a,b}}
\newcommand{\Rabl}{{\mathcal R}_{a,b,\lambda}}
\newcommand{\Grlx}{G_{r,\lambda}(X)}
\newcommand{\Grltau}{G_{r,\lambda}(\tau)}
\newcommand{\Grlrol}{G_{r,\lambda}(\ro_\lambda)}

\newcommand{\asb}{\frac{a}{b}}
\newcommand{\Psilr}{\Psi_{ \lambda}(r)}
\newcommand{\Psil}{\Psi_{ \lambda}}

\newcommand{\unM}{\{1,\ldots,M\}}
\newcommand{\unMmu}{\{1,\ldots,M-1\}}
\newcommand{\Chd}{C}
\newcommand{\Ctihd}{D}

\newcommand{\unppu}{\{1,\ldots,p+1\}}
\newcommand{\calN}{{\mathcal N}}
\newcommand{\calNpr}{{\mathcal N}'}
\newcommand{\Card}{{\rm Card}}

\newcommand{\cstun}{ 6\cdot 10^{-6} \,  }

\newcommand{\dun}{23.0000987}
\newcommand{\dde}{888.376706   \, \,   }
\newcommand{\dtr}{17.068934}
\newcommand{\cun}{0.003261 \, }
\newcommand{\cde}{605.44564  \,  }
\newcommand{\ctr}{603.22318  \, }
\newcommand{\cqu}{3  \cdot 10^{-6}\, }
\newcommand{\deuxcqu}{6 \cdot10^{-6}\, }
\newcommand{\cci}{6 \cdot 10^{-6}  \,  }
\newcommand{\csi}{0.993477 \, }

\title{Distribution of irrational zeta values}
\author{St\'ephane Fischler}
\date{\today}

\begin{document}

\newcommand{\pb}[1]{{\bf #1}}

\maketitle

\begin{abstract}
In this paper we refine Ball-Rivoal's theorem by proving that for any odd integer $a$ sufficiently large in terms of $\eps>0$, there exist 
$[ \frac{(1-\eps)\log a }{1+\log 2}]$ odd integers $s$ between 3 and $a$, with distance at least $a^\eps$ from one another, at which Riemann zeta function takes $\Q$-linearly independent values. As a consequence, if there are very few integers $s$ such that $\zeta(s)$ is irrational, then they are rather evenly distributed.

The proof involves series of hypergeometric type estimated by the saddle point method, and the generalization to vectors of Nesterenko's linear independence criterion.
\end{abstract}


\bigskip


\section{Introduction}

Conjecturally, the values of Riemann zeta function at odd integers $s\geq 3$ are irrational, and together with 1 they are linearly independent over the rationals.  However very few results are known in this direction. After Ap\'ery's breakthrough, namely the proof \cite{Apery} that $\zeta(3)$ is irrational, the next major result is due to Ball-Rivoal  (\cite{BR}, \cite{RivoalCRAS}):

\begin{Th}[Ball-Rivoal] \label{thBR} Let $\eps >0$,  and $a$ be an odd integer sufficiently large with respect to $\eps$. Then the $\Q$-vector space 
\begin{equation} \label{eqthBR}
\Span_\Q(1,\zeta(3),\zeta(5),\ldots,\zeta(a))
\end{equation}
has dimension at least $\frac{1-\eps}{1+\log 2}\log(a )$.
\end{Th}

Except when  $a$ is  bounded, this is the only known linear independence result on the values $\zeta(s)$ for odd $s\leq a$. Trying only to find integers $s$ such that $\zeta(s)$ is irrational, the following result of Zudilin  (Theorem 0.2 of \cite{Zudilincentqc}) has also to be mentioned.

\begin{Th}[Zudilin] \label{thhuit} For any odd integer $d \geq 1$, at least one of the numbers
$$\zeta(d+2), \hspace{0.3cm} \zeta(d+4), \hspace{0.3cm} \zeta(d+6), \hspace{0.3cm} \ldots,   \hspace{0.3cm} \zeta(8d-1)$$
is irrational.
\end{Th}

The purpose of the present paper is to prove results on the distribution of (provably) irrational (or linearly independent) zeta values. For instance, given a large odd integer $a$, Theorems \ref{thBR} and \ref{thhuit} don't exclude the possibility that $1$, $\zeta(3)$, $\zeta(5)$, \ldots, $\zeta(N)$ are $\Q$-linearly independent, with $N = [\frac{\log a}{1+\log 2}]$, and $\zeta( N+2)$,  $\zeta( N+4)$,  \ldots, $\zeta(a)$ are all rational multiples of $\zeta(3)$. More generally, there might exist a few small blocks of consecutive odd integers among which one has to take the integers $s\leq a $ so that the values $\zeta(s)$ make up a basis of the $\Q$-vector space \eqref{eqthBR}. We prove that this cannot happen, as the following result shows.

\begin{Th}\label{th70} Let $\eps >0$,  and $a\geq d \geq 1$ be such that $0 < \eps \leq 1/20$ and $a\geq \eps^{-12/\eps} d$. Then there exist odd integers  $\sigma_1$, \ldots, $\sigma_N$ between $d$ and $a$, with $N = [\frac{1-\eps}{1+\log 2}\log(a/d)]$,  such that:
\begin{itemize}
\item 1, $\zeta(\sigma_1)$, \ldots, $\zeta(\sigma_N)$ are linearly independent over the rationals.
\item  For any $i\neq j$, we have $|\sigma_i - \sigma_j| > d$.
\end{itemize}
\end{Th}

Taking $d = a^\eps$ in this result, we obtain Theorem \ref{thBR} with two additional properties: linearly independent zeta values with distance at least $a^\eps$ from one another, and an explicit value $a(\eps)$ such that the conclusion of Theorem \ref{thBR} holds for any $a\geq a(\eps)$. The latter could have been derived from Ball-Rivoal's proof (\cite{BR}, \cite{RivoalCRAS}), whereas the former is the central new result of the present paper.

\bigskip

Coming back to arbitrary values of $d$, one may weaken the conclusion $|\sigma_i - \sigma_j| > d$ of Theorem \ref{th70} to $\sigma_i  > d $, discarding at most one zeta value $\zeta(\sigma_j)$. This yields the following corollary, in which for simplicity we omit the explicit relations of Theorem \ref{th70} on $\eps$, $d$, $a$.

\begin{Cor} \label{cordiode}
Let $\eps > 0$.  Let $a \geq 3 $ and $d \geq 1$ be odd integers such that $a/d$ is sufficiently large (in terms of $\eps$). Then 
$$\dim_\Q \Span_\Q (\zeta(d), \zeta(d+2), \zeta(d+4), \ldots, \zeta(a)) \geq  \frac{  (1-\eps)\log(a/d)}{1+\log 2}.$$
\end{Cor}

\bigskip

Moving now to bounded values of $a$, Ball and Rivoal have proved (\cite{BR}, \cite{RivoalCRAS}) that \eqref{eqthBR} has dimension at least 3 for $a = 169$. This numerical value has been improved to 145 by Zudilin \cite{Zudilincentqc}, and to 139 in \cite{SFZu}. We obtain the following  result in the spirit of Theorem \ref{thhuit} and Corollary \ref{cordiode}.

\begin{Th} \label{th145}
For any odd integer $d\geq 1$ there exist odd integers $\sigma_1,\sigma_2$ with 
$$d+2 \leq \sigma_1 < \sigma_2 \leq 151\,d, \hspace{1cm} \sigma_2 > \sigma_1 + \cstun d,$$
  such that 1, $\zeta(\sigma_1)$ and $\zeta(\sigma_2)$ are $\Q$-linearly independent.
\end{Th}

This result is new   for any $d \geq 3$,  even if $\sigma_2 > \sigma_1 + \cstun d$  is replaced with $\sigma_2 > \sigma_1 $.  The numerical value 151 (instead of 145 or 139) comes from the fact that some estimates are slightly worse when $d$ is large than for $d=1$.

\bigskip

Now let us move from linear independence to irrationality of zeta values. Ball-Rivoal's theorem yields an increasing sequence $(u_i)_{i\geq 1}$ of odd integers such that $\zeta(u_i)\not\in\Q$ for any $i$, and $\limsup u_i^{1/i}\leq 2e$; for instance it is enough to denote by $u_i$ the $i$-th odd integer $s\geq 3$ such that $\zeta(s)\not\in\Q$. The existence of such a sequence with $\lim  u_i^{1/i}= 2e$ can be deduced from Corollary \ref{cordiode} (by following the proof of Corollary \ref{cor82} below). Actually, using Theorem \ref{th70} we obtain the following result, in which the odd integers $u_i$ are quite distant from one another.

\begin{Cor}\label{cor82} Let $\eps  $ be a positive real number such that $\eps \leq 1/20$; put $\eta = \eps^{15/\eps}$. Then there exists an increasing sequence $(u_i)_{i\geq 1}$ of odd integers, depending only on $\eps$, with the following properties:
\begin{itemize}
\item For any $i\geq 1$, $\zeta(u_i)$ is an irrational number.
\item For any $i\geq 1$, we have $u_{i+1}/u_i > 1+\eta$.
\item For any $i\geq 1$, we have $\eta (2e)^{(1+\eps)i} < u_i < \eta^{-1} (2e)^{(1+\eps)i} $.
\item For any $a \geq \eta^{-1/\eps}$ we have $s_N \leq a$, where $N $ is the integer part of $\frac{1-2\eps}{1+\log 2}\log a$.
\end{itemize}
\end{Cor}

The point here is that the lower bound  $u_{i+1}  > (1+\eta)u_i$ is much stronger than the one of Theorem \ref{th70}, namely $u_{i+1} > u_i+d$ where $d$ has to be comparable to $\log a $ (and therefore to $\log u_i$, at least for most values of $i$) in order to keep a proportion of irrational zeta values as large as in Ball-Rivoal's result. Note that in this respect, Corollary \ref{cor82} refines on Ball-Rivoal's theorem, Theorem \ref{th70} and Corollary \ref{cordiode}, if in these statements the linear independence with 1 is replaced with irrationality.

\bigskip

If we imagine that only $[ \frac{ 1-\eps  }{1+\log 2}\log a]$ odd integers $s \leq a$ are such that $\zeta(s)\not\in\Q$, then (up to a few exceptions) these are the odd integers $u_1,\ldots , u_N$ of Corollary \ref{cor82}; in particular they are rather well distributed.

\bigskip

To conclude this introduction we mention the following result, analogous to the one of \cite{Pilehroodssums} concerning the numbers $\lambda_0 \zeta(s) + \lambda_1 s  \zeta(s+1)$ (see also Th\'eor\`eme 2 of \cite{FR}).

\begin{Th} \label{thzetafacile}
Let $\eps$, $a$, $d$ be as in Theorem \ref{th70}. Let $\lambda_0,\ldots,\lambda_d$ be real numbers, not all zero. Then the real numbers
\begin{equation}\label{eqfac}
\lambda_0\zeta(s) +  \lambda_1  \combi{s+1}{2} \zeta(s+2) + \lambda_2  \combi{s+3}{4} \zeta(s+4) + \ldots + \lambda_d \combi{s+2d-1}{2d} \zeta(s+2d), 
\end{equation}
for odd integers $s$ between $d$ and $a$, span a $\Q$-vector space of dimension at least $[\frac{1-\eps}{1+\log 2}\log(a/d)]$.
\end{Th}

\begin{Cor}\label{corzetafacile}
Let $d\geq 1$ and $\lambda_0,\ldots,\lambda_d$ be real numbers, not all zero. Then the number \eqref{eqfac} is irrational for infinitely many odd integers $s$.
\end{Cor}

\bigskip

The proofs of the results stated in this introduction rely on a classical construction of linear forms in zeta values, namely 
\begin{equation} \label{eqintroln}
  \sum_{k =  1}^\infty  \frac{{\rm d}^{\beta-1}}{{\rm d}t ^{\beta-1}} \Big( \frac{(k- 2r n)_{2rn}^b (k+2n+1)_{2rn}^b }{(k)_{2n+1}^a} \Big) 
  \end{equation}
for suitable parameters $a$, $b$, $r$, $n$ (see \S  \ref{subsec61} for details). These are  linear forms  small at several points, and the generalization to vectors \cite{SFnestsev} of Nesterenko's linear independence criterion \cite{Nesterenkocritere} enables one to deduce a lower bound on  the rank of a family of vectors of which the coordinates involve zeta values; this lower bound is our main Diophantine result, stated as Theorem \ref{thdiogen} in \S \ref{subsecenonce}. We would like to outline three main tools used in implementing this strategy, which seem to be new in this context and may be of independent interest:
\begin{itemize}
\item A general result from linear algebra, namely Proposition \ref{propalglin} in \S \ref{secalglin}, enables one to deduce the linear independence of zeta values well apart from one another from the lower bound of Theorem \ref{thdiogen}. This proposition could be used in any other context where the generalization to vectors \cite{SFnestsev} of Nesterenko's linear independence criterion is applied, since it is completely independent from any specific property of Riemann zeta function.
\item It turns out that for $1\leq\beta\leq b-2$, the asymptotic behavior as $n\to\infty$ of the linear forms \eqref{eqintroln} does not depend on $\beta$. Since the linear independence criterion requires linear forms with pairwise distinct asymptotics, we consider  suitable linear combinations
 of the linear forms \eqref{eqintroln}.
\item When applying the saddle point method, as Zudilin did in this context, we have to check that we don't take the real part of a quantity of which the argument tends to $\pi/2 \bmod \pi$ (otherwise we obtain only an upper bound on the upper limit, which is not sufficient to apply the criterion). This is usually done by numerical computations, but we cannot do it here because the parameters vary. Therefore we allow the parameter $r$ to be a {\em rational} number, and prove that this argument tends to  $\pi/2 \bmod \pi$ only for finitely many values of $r$ (namely zeroes of an analytic function). By right-continuity, the output of this method is the same as if this problem had  never occurred. We believe this trick could be applied to other situations where the same problem arises.
\end{itemize}

\bigskip

The structure of this text is as follows. We recall in \S \ref{seccrit} the linear independence criterion, and state in \S \ref{subsecenonce} our main Diophantine result. Sections \ref{subsec61} to \ref{subsec61bis} are devoted to its proof, starting with a sketch and concluding with the details. At last we deduce in \S \ref{seccsq} the results stated in this introduction, starting in \S \ref{secalglin} with the above-mentioned general  proposition of  linear algebra.

\section{The linear independence criterion} \label{seccrit}

Our results are based upon the following criterion, in which    $\R^p$ is  endowed with its canonical scalar product  and the corresponding norm.

\begin{Th} \label{thcritere}
Let $1\leq k \leq p-1$,  and   $e_1,\ldots,e_k \in \R^p$. 

Let $\tau_1>\ldots>\tau_k >0 $ be   real numbers. 

Let $\om_1,\ldots,\om_k, \pha_1,\ldots,\pha_k$ be  real numbers such that $\pha_j \not\equiv \frac{\pi}2 \bmod \pi$ for any $j$. 

Let $(Q_n)_{n \geq 1}$  be an increasing sequence of positive integers, such that  $Q_{n+1} = Q_n^{1+O(1/n)}$.

For any $n\geq 1$, let $ L_n = \ell_{1,n}X_1 + \ldots + \ell_{p,n}X_p$  be a linear form on $\R^p$, with integer coefficients $\ell_{i,n}$ such  that, as $n\to\infty$:
\begin{equation} \label{eqhypasyth}
|L_n(e_j) | = Q_n^{-\tau_j+o(1)} |\cos(n \om_j + \pha_j) + o(1)|  \mbox{ for any } j \in \unk,
\end{equation}
and 
$$\max_{1 \leq i \leq p} |\ell_{i,n} | \leq Q_n^{1+o(1)}.$$

Let $M \in {\rm Mat}_{k,p}(\R)$ be the matrix of which $e_1,\ldots,e_k \in \R^p$ are the rows; denote by  $C_1,\ldots,C_p \in \R^k$ its columns. 
Then 
$$\rk_\Q ( C_1 ,\ldots, C_p ) \geq k +  \tau_{1} + \tau_{2} + \ldots + \tau_k$$
where $\rk_\Q ( C_1 ,\ldots, C_p )$ is the rank of the family $ ( C_1 ,\ldots, C_p )$ in $\R^k$ seen as a $\Q$-vector space.
\end{Th}

This result is proved in \cite{SFnestsev}, with a little difference: instead of  $\pha_j \not\equiv \frac{\pi}2 \bmod \pi$ for any $j$,  it is assumed that there exist infinitely many integers $n$ such that,  for any $j \in \unk$, $n \om_j + \pha_j \not\equiv \frac{\pi}2 \bmod \pi$. However the former assumption implies the latter. Indeed, let $J$ be the set of all $j \in \unk$ such that $\om_j/\pi\in\Q$. Let $d$ be a common denominator of the numbers $\om_j/\pi$, $j\in J$; if $J = \emptyset$ we let $d=1$. If $n$ is a multiple of  $d$ then $n \om_j + \pha_j  \equiv  \pha_j  \not\equiv \frac{\pi}2 \bmod \pi$ for any $j\in J$. Moreover for each $j\in\unk\setminus J$ there is at most one integer $n$ for which  $n \om_j + \pha_j  \equiv \frac{\pi}2 \bmod \pi$. Therefore there exist infinitely many integers $n$ such that,  for any $j \in \unk$, $n \om_j + \pha_j \not\equiv \frac{\pi}2 \bmod \pi$. 

\bigskip

\begin{Rem} \label{Rk1}
In the proof of Theorem \ref{th145} we shall use the following refinement (see Corollary 1 of \cite{SFnestsev}). Under the assumptions of Theorem \ref{thcritere}, let $\pi : \R^k \to \R^t$ be a surjective $\R$-linear map, with $t\geq 1$. Then 
$$\rk_\Q (\pi(C_1),\ldots,\pi(C_p)) \geq t + \tau_{k+1-t} + \tau_{k+2-t} + \ldots + \tau_k$$
where $\rk_\Q (\pi(C_1),\ldots,\pi(C_p))$ is the rank of the family $(\pi(C_1),\ldots,\pi(C_p))$ in $\R^t$ seen as a $\Q$-vector space.
\end{Rem}

\section{The main Diophantine result} \label{seczeta}

In this section we state (in \S \ref{subsecenonce}) and then prove our main Diophantine result, of which all results stated in the introduction will follow. We sketch the proof in \S \ref{subsec61}, and give details in  \S \ref{subsec61bis}. In the meantime, we recall Zudilin's results on the saddle point method (\S \ref{subsecasydebut}) and make two important steps: the construction of an invertible matrix (\S \ref{subsecmatri}), and the proof that only finitely many values of the parameter $r$ lead to an imaginary main part when applying
 the saddle point method (\S \ref{subsecchoicer}).

\subsection{Statement of the result} \label{subsecenonce}

Eventhough its conclusion is more involved than the ones of the results stated in the introduction, the following theorem is  the real Diophantine output of our proof; we refer to   \cite{Gutnik83}, \cite{Gutnik2003} and \cite{Hessami} for results of the same flavour, but providing (under very strict assumptions) the linear independence of the whole set of vectors under consideration. 

\begin{Th} \label{thdiogen}
Let $a $ and $b  $ be  positive odd integers such that $b$ divides $a$ and $a\geq 9b$. Consider the following vectors   in $\R^{(b+1)/2}$:
$$
v_1=\left(\begin{array}{c} \zeta(3) \\ \combitiny{4}2 \zeta(5) \\   \combitiny{6}4 \zeta(7) \\ \vdots \\  \combitiny{b+1}{b-1} \zeta(b+2)     \end{array} \right), 
v_2=\left(\begin{array}{c} \zeta(5) \\ \combitiny{6}2 \zeta(7) \\   \combitiny{8}4 \zeta(9) \\ \vdots \\  \combitiny{b+3}{b-1} \zeta(b+4)    \end{array} \right), 
\ldots, 
v_{(a-1)/2}=\left(\begin{array}{c} \zeta(a) \\ \combitiny{a+1}2 \zeta(a+2) \\   \combitiny{a+3}4 \zeta(a+4) \\ \vdots \\  \combitiny{a+b-2}{b-1} \zeta(a+b-1)    \end{array} \right) ,
$$
and denote by $(u_1,\ldots,u_{(b+1)/2})$ the canonical basis of $\R^{(b+1)/2}$. Then in $ \R^{(b+1)/2}$ seen as a $\Q$-vector space, the family of vectors $(u_1,u_2, \ldots,u_{(b+1)/2}, v_1,v_2,  \ldots, v_{(a-1)/2})$ has rank
greater than or equal to
\begin{equation} \label{eqminodiogen}
\frac{b+1}2 \sup_{r\in \Iab} \Big( 1 - \frac{\log \alrab}{\log \Qrab}\Big)
\end{equation}
where $\Iab$ is the set of all real numbers $r\geq 1$ such that $\frac92 br \log(4r+3)\leq a$, 
$$\alrab = \frac{e^{2(a+b-1)}2^{ 2b(r+1)}}{r^{2(a-2br)}\{r\}^{4b\{r\}}} \mbox{ and } \Qrab = \frac{e^{2(a+b-1)}2^{ 2(a-2b[r])} (2r+1)^{2b(2r+1)}}{ (2\{r\})^{4b\{r\}}};$$
here $[r]$ and $\{r\}$ denote the integer and fractional parts of $r$, respectively.
\end{Th}

If $r$ is an integer then $\{r\} = 0$; in this case  the factors $\{r\}^{4b\{r\}}$ and $(2\{r\})^{4b\{r\}}$ disappear since they are equal to $1$.

 Since $a\geq 9b$ we have $\Iab\neq\emptyset$. However, if $a < 9b$ then  $a\leq 7b$ so that $\Iab=\emptyset$. 
 Of course  Theorem \ref{thdiogen} is interesting when $\alrab < 1$ for some $r\in\Rab$, but it holds also otherwise.

\begin{Rem}\label{Rk2} 
In the proof of Theorem \ref{th145} we shall use the following refinement of Theorem~\ref{thdiogen}, which comes from Remark \ref{Rk1}. 
Let $\pi : \R^{(b+1)/2} \to \R^t$ be a surjective $\R$-linear map, with $t\in\{1,\ldots,(b+1)/2\}$. Then in $\R^{t}$  seen as a $\Q$-vector space, the family of vectors $ \pi(u_1)$, $\pi(u_2)$, \ldots, $\pi(u_{(b+1)/2})$, $\pi(v_1)$, $\pi(v_2)$, \ldots, $\pi(v_{(a-1)/2})$ has rank 
greater than or equal to
$$t\sup_{r\in \Iab} \Big( 1 - \frac{\log \alrab}{\log \Qrab}\Big).$$
This is specially interesting when $\pi$ is defined over $\Q$, because in this case the $\Q$-vector space spanned by $ \pi(u_1)$, $\pi(u_2)$, \ldots, $\pi(u_{(b+1)/2})$ is $\Q^t$, so that we obtain 
$$\rk_\Q(u'_1,u'_2,\ldots,u'_t,\pi(v_1),\pi(v_2) , \ldots,  \pi(v_{(a-1)/2})) \geq t\sup_{r\in \Iab} \Big( 1 - \frac{\log \alrab}{\log \Qrab}\Big),$$
where $(u'_1,\ldots,u'_t)$ is the canonical basis of $\R^t$.  The most interesting example of this situation is when $\pi$ is projection on the last $t$ coordinates; this is the one used in the proof of Theorem \ref{th145} (see \S \ref{secdem145}). It allows one to get rid of $\zeta(3)$, $\zeta(5)$, \ldots, $\zeta(b+2-2t)$ in the entries of the vectors $v_j$'s.
\end{Rem}

\begin{Rem} If $a/b$ is sufficiently large with respect to some $\eps>0$, then the lower bound \eqref{eqminodiogen} is greater than or equal to $\frac{b+1}2\frac{1-\eps}{1+\log 2} \log(a/b)$ (see the proof of Theorem \ref{th70} in \S \ref{secdem70}). The result deduced in this way from Theorem \ref{thdiogen} is new even when $a\to\infty$ and $b$ is fixed (already when $b=3$). 
\end{Rem}

\subsection{Overview of the proof} \label{subsec61}

In this section we construct the linear forms used in the proof of Theorem \ref{thdiogen}, and summarize their properties. Some of them follow from results in the literature, or can be proved easily; the other ones will be proved below. 

\bigskip

Let $a$, $b$,  $n$ be positive integers, and $r \geq 0$ be a rational number,  such that $a$ and $b$ are odd,  $rn$ is an integer, and $2br < a$. We denote by $\cale = \{1,3,5,\ldots,b\}$ the set of all odd integers $\beta$ between 1 and $b$, and for any $\beta\in\cale$ we let 
\begin{equation} \label{eqdefibn}
\ibn =   \frac{ (2n)!^{a-2b[r]} }{( \beta-1)! (2\{r\}n)!^{2b}}   \sum_{t = n+1}^\infty \frac{{\rm d}^{\beta-1}}{{\rm d}t ^{\beta-1}} \Big( \frac{(t-(2r+1)n)_{2rn}^b (t+n+1)_{2rn}^b }{(t-n)_{2n+1}^a}\Big)
\end{equation}
where the derivative is taken at $t$. As usual we denote by $[r]$ and $\{r\}$ the integer and fractional parts of $r$ respectively, and   Pochhammer's symbol is defined by $(\alpha)_k = \alpha (\alpha+1)\ldots(\alpha+k-1)$.
Letting $k=t-n$, we have obviously
$$
\ibn =  \frac{ (2n)!^{a-2b[r]} }{( \beta-1)! (2\{r\}n)!^{2b}}  \sum_{k =  1}^\infty  \frac{{\rm d}^{\beta-1}}{{\rm d}t ^{\beta-1}} \Big( \frac{(k- 2r n)_{2rn}^b (k+2n+1)_{2rn}^b }{(k)_{2n+1}^a} \Big),$$
 where the sum actually starts at $k = 2rn+1$. It is not difficult  to prove that 
$$\ibn = \ltibn +   \ell_{3,n}   \combi{\beta+1}{\beta-1} \zeta(\beta+2) +  \ell_{5,n}  \combi{\beta+3}{\beta-1}  \zeta(\beta+4) + \ldots + \ell_{a,n}   \combi{\beta+a-2}{\beta-1} \zeta(\beta+a-1)$$
with rational numbers $ \ltibn$ and $\ell_{i,n}$ (for odd integers $\beta$ and $i$ such that $1 \leq \beta \leq b$ and $3 \leq i \leq a$). Moreover $d_{2n}^{a+b-1}$ is a common denominator of these rational numbers, where $d_k$ is the least common multiple of 1, 2, \ldots, $k$. Recall (for ulterior use) that $d_k = e^{k+o(k)}$ as $k \to \infty$, an equivalent form of the Prime Number Theorem.

We shall also   need an upper bound on the coefficients of the linear forms $\ibn$, namely
\begin{equation} \label{eqmajocoeffs}
\max\Big( \max_\beta |\ltibn|, \max_i |\ell_{i,n}|\Big)   \leq \Big[ \frac{2^{2(a-2b[r])} (2r+1)^{2b(2r+1)} }{(2\{r\})^{4b\{r\}} } \Big]^{n+o(n)}
\end{equation}
as $n \to \infty$ with $rn\in\Z$. This can be proved easily along the same lines as   Proposition~3.1 of \cite{Zudilincentqc}, where $r$ is assumed to be a positive integer. The denominator$(2\{r\}n)!^{2b}$ in Eq.~\eqref{eqdefibn} is responsible for the factor  $(2\{r\})^{4b\{r\}}$, since  if $\{r\}\neq 0$ Stirling's formula yields  $(2\{r\}n)!^{1/n}  \sim (2\{r\})^{2\{r\}} (n/e)^{2\{r\}}$. Of course $\{r\}= 0$ if $r$ is an integer, and $(2\{r\})^{4b\{r\}}$ should be understood as 1 in this case;  then we have also $(2\{r\}n)!^{2b}=1$ so that this factor disappears in Eq. \eqref{eqdefibn}.

\bigskip

Theorem \ref{thcritere} almost applies  to this setting (see \S \ref{subsec61bis}); the 
 difficulties come from the asymptotic estimates of the linear forms $\ibn$. To begin with,   $I_{3,n}$, $I_{5,n}$, \ldots, $I_{b,n}$ have essentially the same size as $n\to\infty$ (see  the end of    \S \ref{subsecmatri} below), so that the assumption that $\tau_1$, \ldots, $\tau_k$ are pairwise distinct is not satisfied (unless $b\leq 3$). Indeed, $\ibn$ can be estimated asymptotically in terms of complex integrals $\jln$ (defined just before the statement of Lemma \ref{lemzu}, in \S \ref{subsecasydebut} below).  Following the proof of Lemma 2.5 and Corollary 2.1 of  \cite{Zudilincentqc} (in which only the case where $r\in\Z$ and $\beta=b$ is considered), one obtains
\begin{equation} \label{eqdeficz}
  \frac{  (2\{r\}n)!^{2b}}{  (2n)!^{2b\{r\}}} \ibn = \pi^\beta \itibn \frac{(-1)^n (2\sqrt{\pi n})^{a-2rb}2^b}{n^{a-1}} \Big(1+O(n^{-1})\Big) \mbox{ as } n \to\infty \mbox{ with } rn\in\Z,
  \end{equation}
where
\begin{equation} \label{eqdefic}
\itibn = -2 \sum_{\lambda\in\cale}\clb \, \,   \Re\,  \jln
\end{equation}
and the matrix $\matc$ is defined in Lemma \ref{lemeximat} below; here we multiply $\ibn$ by $ \frac{  (2\{r\}n)!^{2b}}{  (2n)!^{2b\{r\}}}$ so that the normalizing factor in Eq. \eqref{eqdefibn} becomes $\frac{(2n)!^{a-2br}}{(\beta-1)!}$. Zudilin has given, using the saddle point method, a precise asymptotic expression for  $|\Re \,  \jln | $ as $n\to\infty$ (under appropriate assumptions, see Lemma \ref{lemzu} below). This expression depends on $\lambda$, but the previous relations imply   that $\ibn$ has the same order of magnitude for all values of $\beta$ (except $\beta=1$); see   the end of  \S \ref{subsecmatri}  for details. In the notation of   Theorem \ref{thcritere}, all values of $\tau_i$ (except one) would be equal, so that this criterion does not apply. 

To overcome this difficulty, we prove in Lemma \ref{lemeximat} that the matrix $\matc$ is invertible. Denoting by $\matd$ the inverse matrix, we consider the following linear combinations of $I_{1,n}$, \ldots, $I_{b,n}$:
$$\sln = \sum_{\beta\in\cale} \dbl \pi^{-\beta} \ibn$$
for  $\lambda\in\cale$ and    $n\geq 1$ such that $rn\in\Z$. Then we have
\begin{equation} \label{eqslnde}
\sln = \chirab \,   \Re\,  \jln \mbox{ with $\chirab\in\R $ such that } \lim_{ n \to\infty \atop rn\in\Z} | \chirab| ^{1/n} = \frac1{\{r\}^{4b\{r\}}} ,
\end{equation}
since 
\begin{equation} \label{eqreloufact}
 \lim_{ n \to\infty \atop rn\in\Z} \Big[ \frac{ (2n)!^{2b\{r\}}}{(2\{r\}n)!^{2b}}\Big]^{1/n} =  \frac1{\{r\}^{4b\{r\}}}
 \end{equation}
using Stirling's formula. Provided  the saddle point method applies as in  Zudilin's paper, it   turns out that the linear forms $\sln$  have   pairwise distinct asymptotic behaviors; this would allow us to apply  Theorem \ref{thcritere} and conclude the proof.

\bigskip

However another problem arises. Using the saddle point method, $| \sln |  = | \chirab    \,   \Re\,  \jln  | $ can be written as the real part of a  quantity for which a very precise asymptotic estimate is known. However the main part of this estimate might  (for some values of $\lambda$) be an imaginary complex number  for any $n$. In this case, one can only derive  an upper bound for $\limsup |\sln|^{1/n}$, and this is not sufficient to apply  Theorem \ref{thcritere}. To overcome this difficulty, we construct in \S \ref{subsecchoicer} a non-zero analytic function (depending only on $a$ and $b$) which vanishes at all rational numbers $r$ for which this main part is  imaginary for some $\lambda$. This provides a finite set $\Rab$ such that $\lim |\sln|^{1/n}$ exists (and can be computed) as soon as $r\not\in\Rab$, under the mild assumptions that $b$ divides $a$,  $a \geq 5b$ and $1\leq r \leq \frac{a-1}{3b}$. Of course we have no way to control this set $\Rab$: we are not even able (except if some additional assumptions are made on $a$ and $b$) to exclude the case where $\Rab$ contains all  integers $r$ between 1 and $\frac{a-1}{3b}$. However, since $\Rab$ is a finite set and we allow $r$ to be a rational number,  this finite number of exceptions has no influence of the result: a rational $r\not\in\Rab$ can be found in any open interval contained in $[1,\frac{a-1}{3b}]$. 

\bigskip

Finally, assuming that  $r\not\in\Rab$, $r\geq 1$ and $\frac92 br \log(4r+3)\leq a$ we can apply Zudilin's results and obtain
as $n\to\infty$ with $rn\in\Z$:
\begin{equation} \label{eqslnun}
| \sln | =\Big( \frac{  \epsral }{\{r\}^{4b\{r\}}}\Big)^{n+o(n)} |  \cos(n \omral  + \phiral) + o(1)|  \mbox{ for any }\lambda\in\cale,
\end{equation}
with 
\begin{equation} \label{eqineg}
0 < \eps_{ 1} < \eps_{ 3} < \ldots < \eps_{ b} \leq   \frac{2^{2b(r+1)}}{r^{2(a-2br)}} < 1
\mbox{ and }  \phiral \not\equiv \frac{\pi}2 \bmod \pi    \mbox{ for any }\lambda\in\cale.
\end{equation}
This enables us to apply Theorem \ref{thcritere}, and deduce Theorem \ref{thdiogen}.

\subsection{Construction of an invertible matrix} \label{subsecmatri}

In this section we prove that the matrix $\matc$ of Eq. \eqref{eqdefic} (see \S \ref{subsec61}) is invertible.  Recall that $b $ is an odd integer, fixed in this section,  and $\cale = \{1,3,5,\ldots,b\}$. As in \cite{Zudilincentqc} we let 
$$\cot_\beta (z) = \frac{(-1)^{\beta-1}}{(\beta-1)!}   \, \,  \frac{d^{\beta-1}}{dz^{\beta-1}} \cot(z),$$
where $\cot(z) = \frac{\cos z}{\sin z}$ is the cotangent function.

\begin{Lemma} \label{lemeximat} There exists a unique matrix $\matc$ such that
\begin{equation} \label{eqeximat}
\sin^b(z) \cot_\beta(z) = \sum_{\lambda\in\cale} \clb \Big( e^{i\lambda z} +  e^{-i\lambda z}\Big) \mbox{ for any $z$ and any $\beta\in\cale$}.
\end{equation}
Moreover the coefficients $\clb$ are rational numbers, and this matrix is invertible.
 \end{Lemma}

Eq. \eqref{eqdefic} can be proved easily with these numbers $\clb$, by following the proof of Lemma~2.5 and Corollary 2.1 of  \cite{Zudilincentqc} (in which only the case where $r\in\Z$ and $\beta=b$ is considered).

 \bigskip

 \Dem of Lemma \ref{lemeximat}: Lemma 2.2 of  \cite{Zudilincentqc} provides, for any $\beta\in\cale$, a polynomial $V_\beta(X)\in\Q[X]$ of degree at most $\beta$ such that 
 $$\sin^\beta(z) \cot_\beta(z) = V_\beta(\cos z) \mbox{ and } V_\beta(-X) = - V_\beta(X).$$
 We let 
 \begin{equation} \label{eqdefiW}
W_{b,\beta}(X) = (1-X^2)^{(b-\beta)/2} V_\beta(X)
 \end{equation}
 so that 
 $$\sin^b(z) \cot_\beta(z) = W_{b,\beta}(\cos z), \hspace{0.8cm}  \deg W_{b,\beta} \leq b \hspace{0.8cm} \mbox{ and }\hspace{0.8cm} W_{b,\beta}(-X) = -W_{b,\beta}(X).$$
 Letting $X = \frac12 ( Y + Y ^{-1})$, the last two properties yield
 \begin{equation} \label{eqintermedeximat}
  W_{b,\beta}\Big( \frac12 ( Y + Y ^{-1})\Big) =   \sum_{\lambda\in\cale} \clb \Big( Y^\lambda + Y^{-\lambda}\Big)
 \end{equation}
for uniquely defined real numbers $\clb$, which are rational and such that Eq. \eqref{eqeximat} holds.

It remains to prove that the matrix $\matc$ is invertible. If it is not then there exist real numbers $\mu_1$, $\mu_3$, \ldots, $\mu_b$, not all zero, such that $\sum_{\beta\in\cale} \mu_\beta\clb =0$ for any $\lambda\in\cale$. Using Eq. \eqref{eqintermedeximat} this implies $\sum_{\beta\in\cale} \mu_\beta W_{b,\beta}(X)=0$, that is 
$$  \sum_{\beta\in\cale} \mu_\beta (1-X^2)^{(b-\beta)/2} V_\beta(X)=0.$$
 Considering the largest $\beta$ such that $\mu_\beta\neq 0$, this is a contradiction because $V_{\beta}(1) = 1$ (as Eq. (2.4) of  \cite{Zudilincentqc} shows by induction on $\beta$). This concludes the proof of Lemma \ref{lemeximat}.

 \bigskip
 
 Let us conclude this section with a remark (which is not directly used in the proofs). We have $V_1(X)=X$ so that $\deg W_{b,1}(X) = b$ and $c_{b,1}^{(b)}\neq 0$. Using Eqns. \eqref{eqdeficz},  \eqref{eqdefic},  \eqref{eqineg},  \eqref{eqreloufact} and Lemma \ref{lemzu} below, we deduce that (under the assumptions of this lemma) 
 \begin{equation} \label{eqremfin}
\limsup_{n\to\infty} | I_{1,n} |^{1/n} = \frac1{\{r\}^{4b\{r\}}}  \limsup_{n\to\infty} | \Re \, J_{b,n} |^{1/n} =  \eps_{  b} 
 \end{equation}
 so that $|I_{1,n}|$ and $|S_{b,n}|$ have the same asymptotic behavior (in particular Eq. \eqref{eqslnun} holds also for $|I_{1,n}|$). On the other hand, for any odd $\beta\geq 3$ we have $\deg V_\beta=\beta-2$ (see Lemma 2.2 of  \cite{Zudilincentqc}) so that $\deg W_{b,\beta}(X) = b-2$ and $c_{b,\beta}^{(b)}= 0$, $c_{b-2,\beta}^{(b)}\neq 0$. As above, under the assumptions of  Lemma \ref{lemzu} we obtain for $\beta\in\{3,5,\ldots,b\}$:
  $$\limsup_{n\to\infty} | I_{\beta,n} |^{1/n} = \frac1{\{r\}^{4b\{r\}}}  \limsup_{n\to\infty} | \Re \, J_{b-2,n} |^{1/n} =  \eps_{ b-2} .$$
Since this value does not depend on $\beta$, the linear independence criterion does not apply directly to the linear forms corresponding to $I_{\beta,n}$, $\beta\in\cale$ (except if $\beta\leq 3$). This is why the linear combinations $\sln$ have been  introduced in \S \ref{subsec61}.
   
\subsection{Application of the saddle point method}  \label{subsecasydebut}

In this section we recall Zudilin's results  \cite{Zudilincentqc}  based on the saddle point method; we try to use the same notation. The main difference is that Zudilin assumes the parameter $r$ to be an integer, whereas we allow rational values of $r$ (because $\Rab$ may contain all integers $r$, see \S \ref{subsec61}). Unless otherwise stated, the proofs of     \cite{Zudilincentqc} generalize directly to this setting.

\bigskip

Let $a \geq 3$ and $b\geq 1$ be   odd integers,  and $r$ be a positive real number such that $ 3br\leq a$. We assume also that 
\begin{equation} \label{eqhyprelou}
\Big(3+\frac1r\Big)^b < \Big(1+\frac1{2r}\Big)^{a+b}.
\end{equation}
This assumption appears at the bottom of p. 503 of   \cite{Zudilincentqc}. Zudilin proves (p. 504) that it holds if
$  r=1 $ or $ r \geq 2$. 
In the proof of Theorem \ref{thdiogen} we shall use the fact that Eq. \eqref{eqhyprelou} holds for any $r\in (0,\frac{a}{3b}]$  if
$ 5b \leq a$.
Indeed this follows from Zudilin's proof if $ r \geq 2$, and for any $r \in (0,2]$ we have
$$\Big(3+\frac1r\Big)^b < \Big(1+\frac1{2r}\Big)^{6b}\leq \Big(1+\frac1{2r}\Big)^{a+b}$$
since the polynomial $(2X+1)^6-  2^6X^5(3X+1)$ takes only positive values on $(0,2]$.

\bigskip

Let us consider the complex plane with cuts along the rays $(-\infty,1]$ and $[2r+1,+\infty)$; for $\tau \in (2r+1,+\infty)$, we denote by $\tau+i0$ the corresponding point on the upper bank of the cut $[2r+1,+\infty)$. We let for $\tau \in \C \setminus ((-\infty,1] \cup [2r+1,+\infty))$:
\begin{eqnarray} 
 &&\hspace{0.5cm} 
 f(\tau) =  b(\tau+2r+1)\log(\tau+2r+1) + b(2r+1-\tau)\log(2r+1-\tau)   \nonumber \\
&&\hspace{-0.5cm} 
 + (a+b)(\tau-1)\log(\tau-1) - (a+b)(\tau+1)\log(\tau+1)+ 2(a-2br)\log  (2) .\label{eq212}
\end{eqnarray}
In this formula all logarithms are evaluated at positive real numbers if  $\tau$ belongs to   the real interval $(1,2r+1)$, and we choose the determinations so that all of them take real values in this case. 

\bigskip

The complex roots of  the polynomial
\begin{equation} \label{eqdefQ}
Q(X) = (X+2r+1)^b (X-1)^{a+b} -  (X-2r-1)^b (X+1)^{a+b} \in \Q[X] 
\end{equation}
 are localized in Lemma 2.7 of  \cite{Zudilincentqc}. They are all simple; exactly one of them, denoted by $\mu_1$, belongs to the real interval $(2r+1,+\infty)$. There are also exactly $(b-1)/2$ roots in the domain $\Re \,  z > 0$, $\Im \, z > 0$; we denote them by $\ro_1$, $\ro_3$, $\ro_5$, \ldots, $\ro_{b-2}$ with $\Re \, \ro_1 < \ldots < \Re \,  \ro_{b-2}$ since these real parts are pairwise distinct. For convenience we let $\ro_b = \mu_1+i0$, and recall that $\Re\, \ro_{b-2} < \Re\, \ro_b$; we shall also use the fact that $f'(\ro_\lambda)=\lambda i \pi$ for any $\lambda\in\cale$. Of course the polynomial $Q$ and the roots $\ro_\lambda $ (for $\lambda\in\cale$) depend on $a$, $b$,  and $r$  but not on $n$.

  For $\tau\in \C \setminus ((-\infty,1] \cup [2r+1,+\infty))$ we let also
$$f_0(\tau) = f(\tau) - \tau f'(\tau).$$
Since $\log(2r+1-(\tau+i0)) = \log(\tau-(2r+1))-i\pi$ with $\tau - (2r+1) > 0$, we have for $\tau \in (2r+1,+\infty)$:
\begin{eqnarray*} 
&&\hspace{0.8cm} 
 f(\tau+i0)  =  b(\tau+2r+1)\log(\tau+2r+1) + b(2r+1-\tau)\log(\tau -  2r-1  )  \\
&&\hspace{-0.8cm}  
  + (a+b)(\tau-1)\log(\tau-1) - (a+b)(\tau+1)\log(\tau+1)+ 2(a-2br)\log  (2) - b (2r+1-\tau)i\pi.
\end{eqnarray*}
This function of $\tau \in  (2r+1,+\infty)$ is increasing on $(2r+1,\mu_1)$, assumes a maximal value at $\tau = \mu_1$, and is decreasing on $(\mu_1,+\infty)$ (see Eq. (2.34) and Corollary 2.2 of  \cite{Zudilincentqc}). Following the second proof of Lemma 3 in \cite{BR}, we obtain:
\begin{equation} \label{eqdefepsa} 
 \lim_{n \to \infty} \frac1n  \log \Big( \frac{(2\{r\}n)!^{2b}}{(2n)!^{2b\{r\}}}  |I_{1,n}|\Big)   = \Re f(\mu_1+i0) =  \Re f_0(\mu_1+i0)  
 \end{equation}
since $f'(\mu_1+i0) = bi\pi \in i\R$; this estimate will be used below to prove that $\eps_{ b} \leq   \frac{2^{2b(r+1)}}{r^{2(a-2br)}}$ (see Eq. \eqref{eq316bis}).  The main difference with the proof of \cite{BR} is the term $2(a-2br)\log 2$ in \eqref{eq212}, which comes from the fact that the integer denoted here by $n$ is actually $2n$ with the notation of \cite{BR}; this has an effect because of the normalization factor $(2n)!^{a-2br}$ that occurs in $ \frac{(2\{r\}n)!^{2b}}{(2n)!^{2b\{r\}}}   I_{1,n}$.

\bigskip

By applying the saddle point method, Zudilin   proves the following result, where the roots $\ro_1,\ro_3,\ldots,\ro_b$ of $Q$ are defined above, $g$ is  defined   on the cut plane $\C \setminus ((-\infty,1] \cup [2r+1,+\infty))$ by
$$g(\tau) = \frac{(\tau+2r+1)^{b/2}(2r+1-\tau)^{b/2}}{(\tau+1)^{(a+b)/2}(\tau-1)^{(a+b)/2}},$$
and for $\lambda\in\cale$ and $\mu\in\R$ with $1 < \mu < 2r+1$ we let 
$$\jln = \frac1{2i\pi} \int_{\mu-i\infty}^{\mu+i\infty} e^{n(f(\tau)-\lambda i \pi\tau)} g(\tau) d\tau.$$

\begin{Lemma} \label{lemzu}
Assume that  $a \geq 3$  and $b \geq 1$ are odd integers, and $r>0$ is a real  number such that $3br \leq a$  and Eq. \eqref{eqhyprelou} holds. Assume also that
\begin{equation} \label{eq261}
\mu_1 \leq 2r+1+\min\Big(\frac{br(r+1)}{2(a+b)}, \frac{r(r+1)}{3(2r+1)}\Big).
\end{equation}
Let $\lambda\in\cale$. Put
  $$\epsral = \exp \Re  f_0(\ro_\lambda), \hspace{.3cm} \omral = \Im  f_0(\ro_\lambda)  ,  \hspace{.3cm}\mbox{  and  } \hspace{.3cm} \phiral = -\frac12 \arg f''(\ro_\lambda) + \arg g(\ro_\lambda),$$
   and assume that
  \begin{equation} \label{eq255}
\mbox{either }\phiral  \not\equiv \frac{\pi}2 \bmod \pi 
\mbox{ or }\omral  \not\equiv 0 \bmod \pi.
\end{equation}
Then we have, as $n\to\infty$,
$$|\Re \, \jln| =  \epsral^{n+o(n)}  | \cos(n \omral + \phiral) +o(1) | .$$
\end{Lemma}

In this lemma, for $\lambda=b$ we have $\ro_b = \mu_1+i0$ so that   $  f_0(\ro_b) =   f_0(\mu_1+i0) \in\R$ and $\om_b\in\pi\Z$, which is consistant with Eq. \eqref{eqdefepsa} and shows that in Eq. \eqref{eqremfin} both upper limits  are actually limits.

\subsection{Finiteness of the exceptional values of $r$} \label{subsecchoicer}

In this section we prove the following result, which will enable us to choose the parameter $r$ in such a way that Lemma \ref{lemzu} applies and provides an asymptotic estimate for $|\Re \,  \jln |$. We keep the notation of \S \ref{subsecasydebut}.

\begin{Lemma} \label{lemanalytic} Let $a$ and $b$ be positive odd integers such that $b$ divides $a$ and $a \geq 5b$. Then there exists a finite set $\Rab$, depending only on $a$ and $b$, 
with the following property: for any real number $r$ such that $1\leq r \leq \frac{a-1}{3b}$ and $r\not\in\Rab$, we have for any $\lambda\in\cale = \{1,3,\ldots,b\}$:
\begin{equation} \label{eqcclanalytic}
-\frac12 \arg f''(\ro_\lambda)+\arg g(\ro_\lambda)\not\equiv\frac{\pi}2\bmod\pi.
\end{equation}
\end{Lemma}

In the notation of Lemma \ref{lemzu}, the conclusion of Lemma \ref{lemanalytic} is $\phiral\not\equiv\frac{\pi}2\bmod\pi$, so that assumption \eqref{eq255} holds.

\bigskip

\Dem of Lemma \ref{lemanalytic}: Let $a,b\geq 1$ be odd integers such that $b$ divides $a$ and $a \geq 5b$. As noticed
at the beginning of \S \ref{subsecasydebut}, Eq. \eqref{eqhyprelou} holds for any real number $r$ with $0 < r \leq \frac{a}{3b}$, so that Zudilin's results  \cite{Zudilincentqc} recalled in  \S \ref{subsecasydebut} apply. For  any odd integer $\lambda$ such that $-b < \lambda \leq b$  we consider the following polynomial:
$$\Grlx= (X+2r+1)(X-1)^{\asb +1}-e^{\lambda i \pi / b} (2r+1-X)(X+1)^{\asb+1}.$$
Since $b$ is odd and the $b$-th roots of $-1$ are the complex numbers $e^{\lambda i \pi / b} $ for odd integers $\lambda$ such that $-b < \lambda \leq b$, we have the following factorization of the polynomial defined in Eq.~\eqref{eqdefQ}:
$$Q(X) = \prod_{-b < \lambda \leq b  \atop \lambda \mbox{ {\tiny odd }}} \Grlx;$$
accordingly $\Grlx$ divides $Q(X)$ for any $\lambda$. We shall use the fact that for $ \tau\in \C \setminus ((-\infty,1] \cup [2r+1,+\infty))$,
$$\Grltau =0 \mbox{ if, and only if, } \exp(\frac1{b}f'(\tau)) = \exp( \lambda i \pi / b);$$
this follows from the formula
$$f'(\tau) =  b\log(\tau+2r+1)- b\log(2r+1-\tau)+(a+b)\log(\tau-1)- (a+b)\log(\tau+1)$$
(see \cite{Zudilincentqc},  Eq. (2.23)). 

Now let us fix $\lambda\in\cale\setminus\{b\}=\{1,3,5,\ldots,b-2\}$. Then we have $f'(\ro_\lambda)  = \lambda i \pi$  so that $\ro_\lambda $ is a root of $\Grlx$. Moreover $\ro_\lambda$ is the only root of $\Grlx$ with a positive real part (see Lemma 2.7 of   \cite{Zudilincentqc}), and $\Grlx$ has only simple roots (because this property holds for $Q(X)$). Therefore when $\lambda\in\cale\setminus\{b\}$ is fixed, $\ro_\lambda$ is an algebraic function of $r$ with no branch point in $(0, \frac{a}{3b}]$: it is a real-analytic function of $r$ on  $(0, \frac{a}{3b}]$. Eventhough we consider  $\ro_\lambda$  as a  function of $r$, we shall continue (for simplicity) to omit this dependence in the notation.

Now we let 
$$\Psilr = \frac{g(\ro_\lambda )^2}{|g(\ro_\lambda )|^2}\frac{|f''( \ro_\lambda  )|}{ f''( \ro_\lambda  ) }$$
for any real number $r$ such that  $0 < r < \frac{a}{3b}$. This function is well-defined because $ g(\ro_\lambda ) $ and $f''( \ro_\lambda  )$ are non-zero (see \cite{Zudilincentqc}, p. 512),  and it is real-analytic on the real interval $(0,\frac{a}{3b})$. Let us compute its limit as $r\to 0$.
 
We have $|   \ro_\lambda -2r-1| \leq |\mu_1-2r-1|< 2r$ using assumption \eqref{eqhyprelou} (see  \cite{Zudilincentqc}, p. 503), so that $\lim_{r\to 0} \ro_\lambda=1$. Let us write $u=o(v)$ whenever $u$ and $v$ are functions of $r$ such that $\lim_{r\to 0}  u/v = 0$, and $u\sim v$ when $u = v+o(v)$; here the parameters $a$, $b$  and $\lambda$ are fixed. Then we have $\ro_\lambda+1 \sim 2$ and $\ro_\lambda+2r+1 \sim 2$, so that taking equivalents in the relation $\Grlrol=0$ yields
\begin{equation} \label{eqintermedasyrho}
(\ro_\lambda-1)^{\asb+1}\sim 2^{\asb}e^{\lambda i \pi / b} (2r+1-\ro_\lambda).
   \end{equation}
This implies $2r-(  \ro_\lambda-1) = o(\ro_\lambda-1)$ so that $\ro_\lambda-1\sim 2r$. Plugging this equivalence into Eq. \eqref{eqintermedasyrho} yields
$$2r+1-   \ro_\lambda \sim 2 e^{-\lambda i \pi / b} r^{\asb+1}$$
so that
$$\ro_\lambda = 1+2r-  2 e^{-\lambda i \pi / b} r^{\asb+1} + o( r^{\asb+1} ).$$
This enables us to compute the following limit as $r\to 0$:
$$g(\ro_\lambda )^2 = \frac{(\ro_\lambda +2r+1)^b(2r+1-\ro_\lambda )^b}{(\ro_\lambda +1)^{a+b}(\ro_\lambda -1)^{a+b}}
\sim \frac{2^{2b}e^{-\lambda i \pi } r^{a+b}}{2^{2(a+b)} r^{a+b}} = -2^{-2a}$$
since $\lambda$ is an odd integer, so that $\lim_{r\to 0}  \frac{g(\ro_\lambda )^2}{|g(\ro_\lambda )|^2} = -1$. In the same way the quantities $(\ro_\lambda -1)^{-1}$, $(\ro_\lambda+1)^{-1}$ and  $(\ro_\lambda +2r+1)^{-1}$ can all be written as $o((\ro_\lambda -2r-1)^{-1})$; since  
\begin{equation} \label{eqfsec}
f''(\ro_\lambda ) = \frac{b}{\ro_\lambda  + 2r+1} + \frac{b}{2r+1-\ro_\lambda } + \frac{a+b}{\ro_\lambda -1}-\frac{a+b}{\ro_\lambda +1}
\end{equation}
(see at the end of the proof of Lemma 2.9 of \cite{Zudilincentqc}),  we have 
$$f''(\ro_\lambda) \sim  \frac{b }{ 2r+1 -  \ro_\lambda} \sim \frac{b}2 e^{\lambda i \pi / b}  r^{-\asb-1} $$ 
so that  $\lim_{r\to 0}  \frac{|f''( \ro_\lambda  )|}{ f''( \ro_\lambda  ) } = e^{-\lambda i \pi / b}$. Finally we have
\begin{equation} \label{eqlimpsi}
\lim_{r\to 0} \Psilr = - e^{-\lambda i \pi / b} \neq -1.
\end{equation}
Now for any $\lambda\in\cale\setminus\{b\}$ we let $\Rabl$ be the set of all $r\in[1,\frac{a-1}{3b}]$ such that $\Psilr =-1$. If $\Rabl$ is infinite for some $\lambda$ then $\Psil+1$ has non-isolated zeros in the segment $[1,\frac{a-1}{3b}]$: this analytic function of $r$ is identically zero on the real interval $(0, \frac{a}{3b})$. This implies $\lim_{r\to 0} \Psilr = - 1$, in contradiction with Eq. \eqref{eqlimpsi}. Therefore $\Rabl$ is a finite set. Let $\Rab$ be the union of these finite sets, as $\lambda$ ranges through $\cale\setminus\{b\}$. Then for any $\lambda\in\cale\setminus\{b\}$ and any $r\in [1,\frac{a-1}{3b}] \setminus\Rab$ we have $\Psilr \neq -1$ so that $ \arg g(\ro_\lambda)^2 - \arg f''(\ro_\lambda)   \not\equiv \pi  \bmod 2 \pi$. Eq. \eqref{eqcclanalytic} follows for $\lambda\in\cale\setminus\{b\}$. For  $\lambda=b$  we have  $\ro_b = \mu_1+i0$ so that $\arg g(\ro_b) = \frac{-b\pi}2$ and $f''(\ro_b) $ is a negative real number (see Eq. \eqref{eqfsec}). Therefore $\pha_b \equiv 0 \bmod\pi$ since $b$ is odd: Eq. \eqref{eqcclanalytic} holds also for $\lambda=b$.  This concludes the proof of Lemma \ref{lemanalytic}.

\subsection{End of the proof of  Theorem \ref{thdiogen}} \label{subsec61bis}

In this section we complete the proof of     Theorem \ref{thdiogen}  and   Remark \ref{Rk2}.  Let $a\geq 3$ and $b \geq 1$ be odd integers such that $b$ divides $a$ and $a\geq 9b$. 
We put $k = (b+1)/2$  and consider the $k$  following vectors in $\R^{(a+b)/2}$:
$$
\begin{array}{rlcccl}
e_1 =& (1,0,0,\ldots,0,&  \zeta(3) , & \zeta(5) ,& \ldots,&  \zeta(a))\\ 
e_2 =& (0,1,0,\ldots,0,  &\combitiny{4}2 \zeta(5)  ,  & \combitiny{6}2 \zeta(7)  , &\ldots, &  \combitiny{a+1}2   \zeta(a+2))\\ 
e_3 =& (0,0,1,\ldots,0,  & \combitiny{6}4 \zeta(7) , & \combitiny{8}4 \zeta(9), & \ldots,   &\combitiny{a+3}4 \zeta(a+4) ) \\
\vdots \\
e_k = & (0,0,0,\ldots,1, &    \combitiny{b+1}{b-1} \zeta(b+2) , & \combitiny{b+3}{b-1} \zeta(b+4), &  \ldots,  &\combitiny{a+b-2}{b-1} \zeta(a+b-1) ).
\end{array}$$
Let $r$ be a rational number such that $r\geq 1$,  $\frac92 br \log(4r+3)\leq a$,  and $r\not\in\Rab$ (where $\Rab$ is the finite set constructed in Lemma \ref{lemanalytic}). Such a rational number exists since $a \geq 9b$  and $\Rab$ is finite. We keep the notation of  \S \ref{subsec61}, and use the results recalled there.

We denote by $X_1$, $X_3$, \ldots, $X_b$, $Y_3$, $Y_5$, \ldots, $Y_a$ the coordinates on $\R^{(a+b)/2}$ and consider, for any $n\geq 1$ such that $rn\in\Z$, the linear form
$$L_n =d_{2n}^{a+b-1}\Big( \lti_{1,n} X_1  + \lti_{3,n}  X_3 +  \ldots + \lti_{b,n}  X_b + \ell_{3,n} Y_3 + \ell_{5,n} Y_5 + \ldots + \ell_{a,n} Y_a\Big) $$
so that  
$$  L_n(e_j) = d_{2n}^{a+b-1} I_{2j-1,n} \mbox{ for any }j\in\unk,$$
that is $L_n(e_{(\beta+1)/2}) = d_{2n}^{a+b-1}  \ibn$ for any $\beta\in\cale = \{1,3,\ldots,b\}$. Using the matrix $\matd$ 
(which is the inverse of $\matc$, see \S\S \ref{subsec61} and \ref{subsecmatri}),  we let
$$e'_{(\lambda+1)/2}  = \sum_{\beta\in\cale} \dbl \pi^{-\beta}  e_{(\beta+1)/2} \in \R^{(a+b)/2} \mbox{ for any } \lambda\in\cale , $$
so that  
\begin{equation} \label{eqliensln}
 L_n(e'_{(\lambda+1)/2 })  = d_{2n}^{a+b-1}  \sln \mbox{ for any } \lambda\in\cale
 \end{equation} 
by definition of $\sln$ (see   \S \ref{subsec61}). To obtain an asymptotic estimate for $ L_n(e'_{(\lambda+1)/2 }) $, it is enough  (using  Eq. \eqref{eqslnde}) to apply Lemma \ref{lemzu}. Let us check the assumptions of this lemma, starting with Eq.  \eqref{eq261}.

Since the map $x\mapsto \frac{\log x}{x} $ is decreasing on the interval $[e, +\infty)$ and $\frac{a+b}{br} \geq \frac{a}{br} \geq \frac92 \log 7 > e$, we have $-\log(\frac{br}{a+b}) = \log( \frac{a+b}{br} ) \leq \frac{2\log(\frac92 \log 7)}{9\log 7} \frac{a+b}{br}$. On the other hand, since $\frac{b}{2(a+b)}\leq \frac1{20}$ we have
$$\log\Big(1+\frac{b}{2(a+b)}\Big) - \log\Big(1+\frac1r+ \frac{b}{2(a+b)}\Big)
\leq \frac{-1}{r} \frac{1}{1+\frac1r+ \frac{b}{2(a+b)}} \leq\frac{-20}{41r}.$$
Therefore we have
\begin{eqnarray*}
&&\hspace{-0.7cm}
b\log\Big(4r+2+\frac{br}{ a+b }\Big) + (a+b) \Big[ \log\Big(1+\frac{b}{2(a+b)}\Big) - \log\Big(1+\frac1r+ \frac{b}{2(a+b)}\Big)\Big] - b \log\Big(\frac{br}{a+b}\Big)\\
& &\hspace{0.7cm}
\leq  \frac{a+b}r\Big( \frac29 - \frac{20}{41} + \frac{2\log(\frac92 \log 7)}{9\log 7} \Big)<0
\end{eqnarray*}
so that
$$Q\Big(2r+1+\frac{br}{a+b}\Big)
=  (2r)^{a+b} \Big[ \Big(4r+2+\frac{br}{a+b} \Big)^b \Big(1 + \frac{ b}{2(a+b)}\Big)^{a+b} - \Big( \frac{br}{a+b}\Big)^b   \Big(1 + \frac1{r} + \frac{b }{2(a+b)}\Big)^{a+b}\Big]  < 0,
$$
where $Q$ is the polynomial defined in Eq. \eqref{eqdefQ}. Since $Q(2r+1)>0$ and   $\mu_1$ is the only root of $Q$  in the real interval $(2r+1,+\infty)$, we obtain
$$\mu_1 < 2r+1+\frac{br}{a+b} \leq 2r+1+ \frac{br(r+1)}{2(a+b)} \leq 2r+1 +   \frac{r(r+1)}{3(2r+1)}  $$
since $6br\leq a$. Therefore  Eq.  \eqref{eq261} holds. 

\bigskip

Moreover  Lemma \ref{lemanalytic}  yields $\phiral\not\equiv\frac{\pi}2\bmod\pi$, so that  assumption \eqref{eq255} holds. As noticed at the beginning of \S \ref{subsecasydebut}, Eq. \eqref{eqhyprelou} holds since $a\geq 5b$. 
Therefore Lemma \ref{lemzu} applies, and provides real numbers  $\epsral$, $ \omral$ and $  \phiral$. The asymptotic estimate \eqref{eqslnun} is an immediate consequence of Eq. \eqref{eqslnde}. The inequalities 
\begin{equation} \label{eq316bis}
0 < \eps_{ 1} < \eps_{ 3} < \ldots < \eps_{ b}  \leq   \frac{2^{2b(r+1)}}{r^{2(a-2br)}}
\end{equation} 
are a consequence of   Lemma 2.10 of  \cite{Zudilincentqc}, except for the last one that we prove now,  following the second proof of Lemme 3
of \cite{BR}.  Since $k+(2r+1)n < 2^{1+1/r}k$  for any $k > 2rn$,  we have  
\begin{eqnarray*}
 (2n)!^{a-2br}  \frac{(k-2rn)_{2rn}^b (k+2n+1)_{2rn}^b }{(k+1)_{2n}^a} 
 &<& (2n)^{2n(a-2br)} \frac{k^{2brn} (2^{1+1/r}k)^{2brn}}{k^{2an}} \\
 &=& \Big( \frac{2n}{k}\Big)^{2n(a-2br)} 2^{2brn(1+1/r)} < \Big[  \frac{2^{2b(r+1)}}{r^{2(a-2br)}} \Big]^n.
\end{eqnarray*}
This yields
$$ \frac{  (2\{r\}n)!^{2b}}{(2n)!^{2b\{r\}}} |I_{1,n}| \leq \Big[  \frac{2^{2b(r+1)}}{r^{2(a-2br)}} \Big]^n \sum_{k=2rn+1}^{+\infty} \frac1{k^a}$$
so that  $ \eps_{ b} \leq   \frac{2^{2b(r+1)}}{r^{2(a-2br)}}$ (using Eq. \eqref{eqdefepsa} and the fact that $\eps_{ b} =\exp  \Re\, f_0(\mu_1+i0)$). 

\bigskip

We are now in position to apply the linear independence criterion, namely~Theorem \ref{thcritere}. We let 
$$\tau_{(\lambda+1)/2} = \frac{-\log(e^{2(a+b-1)}\epsral  \{r\} ^{-4b\{r\}})}{\log \Qrab}$$ 
 for any $ \lambda\in\cale$, and Eq. \eqref{eq316bis} yields 
 \begin{equation} \label{eqtaudiff}
 \tau_1>\tau_2>\ldots>\tau_k \geq - \frac{\log \alrab}{\log \Qrab}.
 \end{equation} 
 Now  let $\delta$ denote the denominator of the rational number $r$, and $Q_n = \Qrab^{\delta n}$ for any $n\geq 1$. Then Theorem \ref{thcritere} applies to the linear forms $L_{\delta n }$, $n\geq 1$, at the points $e'_1,\ldots, e'_k$,   using (among others) Eqns. \eqref{eqmajocoeffs}, \eqref{eqliensln} and \eqref{eqslnde}, and Lemmas \ref{lemzu} and \ref{lemanalytic}. The columns $C_1,\ldots,C_{(a+b)/2}$  are exactly the vectors denoted by $u_1,\ldots,u_{(b+1)/2},v_1,\ldots,v_{(a-1)/2}$ in  Theorem \ref{thdiogen}. Using Eq. \eqref{eqtaudiff} we obtain in this way the lower bound \eqref{eqminodiogen} for $r\not\in\Rab$; since $\Rab$ is a finite set, the supremum is the same by right-continuity. This concludes the proof of  Theorem \ref{thdiogen}; Remark \ref{Rk2} can be proved in the same way, using   Remark \ref{Rk1} stated after Theorem \ref{thcritere}.

\section{Proof of the Diophantine consequences} \label{seccsq}

In this section we deduce from Theorem \ref{thdiogen} all results stated in the introduction; the main tool is a result coming from linear algebra, stated and proved  in \S \ref{secalglin}.

\subsection{A  linear algebra result} \label{secalglin}

We state in this section one of the main tools in the proof of the results stated in the introduction. It enables one to deduce from a lower bound on the rank of a family of vectors $(v_1,\ldots,v_N)$, such as the one provided by Theorem \ref{thdiogen}, the existence of linearly independent entries of the vectors $v_j$ which are not too close from one another. We state it in a general form, dealing with any vector space $E$ on a field $\K$. We shall apply it with $E = \R/\Q$ and $\K=\Q$: real numbers have linearly independent images in $\R/\Q$ if, and only if,  together with 1 they are  $\Q$-linearly independent in $\R$. We hope this result can be used in other contexts (not involving Riemann zeta function), to take advantage of the lower bound provided by the linear independence criterion.

\bigskip
 
 To state the result, we fix $k\geq 1$, $N \geq 1$, and we let $[\lambda_{i,j}]_{1\leq i \leq k, 1 \leq j \leq N}$ be a $k\times N$ matrix with entries in $\K$ and $\xi : \llbracket 1,N+k-1\rrbracket \to E$ be a map. We consider the following vectors in the $\K$-vector space $E^k$:
$$
v_1 = \left( \begin{array}{c}  \lambda_{1,1}\xi(1) \\ \lambda_{2,1}\xi(2) \\ \vdots \\ \lambda_{k,1}\xi(k)\end{array}\right), 
\hspace{0.8cm}
v_2 = \left( \begin{array}{c}\lambda_{1,2}\xi(2) \\ \lambda_{2,2}\xi(3) \\ \vdots \\ \lambda_{k,2}\xi(k+1)\end{array}\right), 
\hspace{0.6cm}
\ldots, 
\hspace{0.6cm}
v_N = \left( \begin{array}{c}\lambda_{1,N}\xi(N) \\ \lambda_{2,N}\xi(N+1) \\ \vdots \\ \lambda_{k,N}\xi(N+k-1)\end{array}\right).
$$

\bigskip

\begin{Prop}\label{propalglin} Let $\delta \geq 0$ and $p,q\geq 0$ be such that
$$\rk_\K(v_1,\ldots,v_N) > (k+4\delta)(p+q-1).$$
Then for any $ m_1,\ldots,m_q\in \llbracket  1, N+k-1\rrbracket $ there exist $n_1,\ldots,n_p \in  \llbracket  1, N+k-1\rrbracket $ with the following properties:
\begin{itemize}
\item $\xi(n_1)$, \ldots, $\xi(n_p)$ are $\K$-linearly independent.
\item For any $i,j\in\unp$ with $i\neq j$, $|n_i-n_j| > \delta$.
\item For any $i \in \unp$ and any $j\in\unq$, $|n_i-m_j| > \delta$.
\end{itemize}
\end{Prop}

\bigskip

The integer $q$ plays a crucial role in the proof of this proposition, but in this paper we apply it only with $q=0$.

\bigskip

With $\delta=0$, Proposition \ref{propalglin} can be proved easily. Indeed, let $L$ denote the set of indices $\ell$ such that some $\xi(m_j)$, $1\leq j \leq q$, appears in an entry of $v_\ell$; then $L = \llbracket 1,N \rrbracket \cap \cup_{j=1}^q \llbracket  m_j-k+1, m_j\rrbracket$ so that $\Card \, L \leq kq$. Therefore the family $(v_\ell)_{\ell\not\in L}$ has rank greater than $k(p-1)$. Now letting $F$ denote the $\K$-subspace of $E$ generated by the numbers $\xi(n)$ for $n\in  \llbracket 1,N+k-1 \rrbracket\setminus\{m_1,\ldots,m_q\}$, we have $v_\ell \in F^k$ for any $\ell\not\in L$ so that $\dim(F^k) > k(p-1)$ and $\dim F \geq p$. This concludes the proof of Proposition \ref{propalglin} in this case.

\bigskip

To prove Proposition \ref{propalglin} when $\delta>0$, we apply $p$ times the following result, with $R_j = \max(1,m_j-\delta)$ and $S_j = \min(N+k-1, m_j+\delta)$. 

\begin{Lemma} \label{lemalglin}
Let $\delta \geq 1$ and $p,q\geq 0$ be such that
$$\rk_\K(v_1,\ldots,v_N) > (k+4\delta)(p+q).$$
Let $R_1,\ldots,R_q,S_1,\ldots,S_q \in  \llbracket  1, N+k-1\rrbracket $ be such that $R_j \leq S_j \leq R_j+4\delta$ for any $j\in\unq$, and put
$$\calN =  \llbracket  1, N+k-1\rrbracket  \setminus \bigcup_{j=1}^q \llbracket R_j,S_j\rrbracket .$$
Let $n_1,\ldots,n_p\in\calN$ be such that $\xi(n_1)$, \ldots, $\xi(n_p)$ are $\K$-linearly independent and $|n_i-n_j| > \delta$ for  any $i,j\in\unp$ with $i\neq j$. Then there  exist $n'_1,\ldots,n'_{p+1} \in \calN$ such that:
\begin{itemize}
\item $\xi(n'_1)$, \ldots, $\xi(n'_{p+1})$ are $\K$-linearly independent.
\item For any $i,j\in\unppu$ with $i\neq j$, $|n'_i-n'_j| > \delta$.
\item $\Span_\K(  \xi(n_1) , \ldots,  \xi(n_p) ) \subset \Span_\K( \xi(n'_1) , \ldots,  \xi(n'_{p+1}) ).$
\end{itemize}
\end{Lemma}

\Dem of Lemma \ref{lemalglin}: We let 
$$\calNpr = \llbracket 1,N+k-1\rrbracket  \setminus \bigcup_{j=1}^q \llbracket R_j-k+1,S_j\rrbracket $$
and argue by induction on $p$. If $p=0$, the assumption $\rk (v_1,\ldots,v_N) > (k+4\delta) q $ yields $\rk\{v_n, n \in \llbracket 1,N\rrbracket \cap \calNpr\}> 0$ since 
\begin{equation} \label{eqmajoNpr}
\Card(\llbracket 1,N\rrbracket \cap\calNpr)\leq \sum_{j=1}^q (S_j-R_j+k)\leq (k+4\delta)q.
\end{equation}
Therefore $v_n\neq 0$ for some $n \in \llbracket 1,N\rrbracket \cap \calNpr$; there exists $n'_1\in\llbracket n,n+k-1\rrbracket \subset\calN$ such that $\xi(n'_1)\neq 0$. This concludes the proof of Lemma \ref{lemalglin} if $p=0$.

Assume this lemma holds for any $p'\leq p-1$, with $p\geq 1$, and let us prove it for $p$. Consider the vector subspace $F$ of $E$ generated by the elements $\xi(n)$, for $n\in\calN$ such that $|n-n_i|>\delta$ for any $i\in\unp$. If $F$ is not contained in 
 $\Span (  \xi(n_1) , \ldots,  \xi(n_p) )$, we take $n'_1=n_1$, \ldots, $n'_p = n_p$ and there exists $n'_{p+1} \in\calN$ such that $|n'_{p+1}-n_i| > \delta$ for any $i \in\unp$ and $\xi(n'_{p+1})\not\in\Span (  \xi(n_1) , \ldots,  \xi(n_p) )$; the lemma follows at once in this case. Therefore we assume from now on that $F\subset \Span (  \xi(n_1) , \ldots,  \xi(n_p) )$.
 
 Now we have $\rk(v_1,\ldots,v_N) >  (k+4\delta)(p+q)$ so that Eq. \eqref{eqmajoNpr} yields
 $$\rk\{v_n, n \in \llbracket 1,N\rrbracket \cap \calNpr\}> (k+4\delta) p > kp = \dim_\K\Big( \Span (  \xi(n_1) , \ldots,  \xi(n_p) )\Big)^k.$$
Therefore $v_n\not\in \Big( \Span (  \xi(n_1) , \ldots,  \xi(n_p) )\Big)^k$  for some $n \in \llbracket 1,N\rrbracket \cap \calNpr$: there exists $s\in\llbracket n,n+k-1\rrbracket \subset\calN$ such that $\xi(s)\not\in  \Span (  \xi(n_1) , \ldots,  \xi(n_p) )$. Since $F\subset \Span (  \xi(n_1) , \ldots,  \xi(n_p) )$ we have $\xi(s)\not\in F$ so that $|s-n_i| \leq \delta$ for some $i\in\unp$, by definition of $F$. Since $n_1$, \ldots, $n_p$ play symmetric roles we may assume that $i=1$. Let us distinguish between two cases.

\bigskip

$\bullet$ To begin with, let us consider the case where $|s-n_i| > \delta$ for any $i\in\{2,\ldots,p\}$; in particular this holds if $p=1$. Then we let $R_{q+1} = \min(s,n_1)-\delta$ and   $S_{q+1} = \max(s,n_1)+\delta$ so that $n_2,\ldots,n_p\not\in\llbracket R_{q+1},S_{q+1}\rrbracket $. Therefore Lemma \ref{lemalglin} applies with $R_1$, \ldots, $R_{q+1}$, $S_1$, \ldots, $S_{q+1}$, and $n_2$, \ldots, $n_p$. This provides integers $n'_2$, \ldots, $n'_p$, $n'_{p+1}$ such that:
\begin{itemize}
\item[$(a)$] $n'_2 , \ldots,  n'_p ,  n'_{p+1} \in \calN\setminus \llbracket R_{q+1},S_{q+1}\rrbracket $,
\item[$(b)$] $\xi(n'_2),\ldots,\xi(n'_p),\xi(n'_{p+1})$ are $\K$-linearly independent,
\item[$(c)$] For any $i,j\in\{2,\ldots,p+1\}$ with $i\neq j$, $|n'_i-n'_j|>\delta$,
\item[$(d)$] $\Span (  \xi(n_2) , \ldots,  \xi(n_p) )\subset \Span (  \xi(n'_2) , \ldots,  \xi(n'_p),\xi(n'_{p+1} ))$.
\end{itemize}
Now $ \xi(n'_2) , \ldots,  \xi(n'_p),\xi(n'_{p+1} )$ are $p$ linearly independent vectors thanks to $(b)$, and  $\xi(n_1)$, \ldots,  $\xi(n_p)$, $\xi(s)$ are $p+1$  linearly independent vectors by construction of $s$. Therefore one can find $n'_1\in\{n_1,\ldots,n_p,s\}$ such that $ \xi(n'_2)$, \ldots,  $\xi(n'_p)$, $\xi(n'_{p+1} )$, $\xi(n'_1)$ are  $p+1$  linearly independent vectors. Assertion $(d)$ above implies $n'_1\not\in\{n_2,\ldots,n_p\}$, so that $n'_1=n_1$ or $n'_1=s$; if possible we choose $n'_1=n_1$. Let us check the conclusions of Lemma  \ref{lemalglin} with $n'_1$,\ldots, $n'_{p+1}$.

Assertion $(a)$ and the construction of $s$ yields $n'_1 ,\ldots, n'_{p+1}\in\calN$; and $ \xi(n'_1)$, \ldots,  $\xi(n'_{p+1})$ are  $\K$-linearly independent by definition of $n'_1$. Given $i,j\in\unppu$ with $i\neq j$, we have $|n'_i-n'_j|>\delta$: this follows from assertion $(c)$ if $i,j\geq 2$, and from $(a)$ if $i=1$ or $j=1$ (by definition of $R_{q+1}$ and $S_{q+1}$, since $n'_1\in\{n_1,s\}$). At last, assertion $(d)$ yields
$$\Span (  \xi(n_1) , \ldots,  \xi(n_p) )\subset\Span (  \xi(n_1) ,\xi(n'_2), \ldots,  \xi(n'_p) , \xi(n'_{p+1}) ).$$
This concludes the proof  of Lemma  \ref{lemalglin}  if $n'_1=n_1$. Otherwise, namely if $n'_1=s$, we have assumed that choosing $n'_1=n_1$ was not possible so that $ \xi(n'_2)$, \ldots,  $\xi(n'_{p+1})$, $\xi(n_1)$ are   linearly  dependent. Using assertion $(b)$ this implies $\xi(n_1)\in \Span (  \xi(n'_2), \ldots,    \xi(n'_{p+1}) )$ so that assertion $(d)$ yields
$$\Span (  \xi(n_1) , \ldots,  \xi(n_p) )\subset \Span (  \xi(n'_2), \ldots,    \xi(n'_{p+1}) ) \subset \Span (  \xi(n'_1), \xi(n'_2),  \ldots,    \xi(n'_{p+1}) ).$$
This concludes the proof  of Lemma  \ref{lemalglin} in the first case.

\bigskip

$\bullet$ Let us move now to the second case: assume there exists $i\in\{2,\ldots,p\}$ such that  $|s-n_i| \leq \delta$. We may assume that $i=2$ has this property. Exchanging $n_1$ and $n_2$ if necessary, we may also assume that $n_1< n_2$; since $n_2-n_1>\delta$ this implies $n_1< s< n_2$. We let $R_{q+1} = n_1-\delta$ and $S_{q+1} = n_2+\delta$, so that $n_3,\ldots,n_p\not\in\llbracket R_{q+1},S_{q+1}\rrbracket $ (because $n_2-n_1\leq |s-n_1|+|s-n_2|\leq 2 \delta$ so that no integer $n\in \llbracket R_{q+1},S_{q+1}\rrbracket $ satisfies both $|n-n_1|> \delta$ and $|n-n_2|> \delta$). Therefore Lemma  \ref{lemalglin} applies with  $R_1$, \ldots, $R_{q+1}$, $S_1$, \ldots, $S_{q+1}$, and $n_3$, \ldots, $n_p$. It provides integers $n'_3 , \ldots, n'_{p+1} \in\calN \setminus \llbracket R_{q+1},S_{q+1}\rrbracket $, and we apply it again with $R_1$, \ldots, $R_{q+1}$, $S_1$, \ldots, $S_{q+1}$, and $n'_3$, \ldots, $n'_{p+1}$. We obtain in this way  integers $n''_3 , \ldots, n''_{p+1}, n''_{p+2} $ such that:
\begin{itemize}
\item[$(a)$] $n''_3 , \ldots,  n''_{p+2} \in \calN\setminus \llbracket  n_1-\delta, n_2+\delta\rrbracket $,
\item[$(b)$] $\xi(n''_3),\ldots, \xi(n''_{p+2})$ are $\K$-linearly independent,
\item[$(c)$] For any $i,j\in\{3,\ldots,p+2\}$ with $i\neq j$, $|n''_i-n''_j|>\delta$,
\item[$(d)$] $\Span (  \xi(n_3) , \ldots,  \xi(n_p) )\subset \Span (  \xi(n'_3) , \ldots, \xi(n'_{p+1} ))   \subset \Span (  \xi(n''_3) , \ldots, \xi(n''_{p+2} ))$.
\end{itemize}
Of course the corresponding  properties  hold  also for $n'_3$, \ldots, $n'_{p+1}$. Now let  us distinguish three cases according to the value of 
$$d = \dim  \Big( \Span  (  \xi(n_1)  ,  \xi(n_2) )  \cap    \Span  (  \xi(n''_3)  , \ldots ,  \xi(n''_{p+2}) )\Big) \in \{0,1,2\}.$$ 
If $d =0$ then we have also
$$ \Span  (  \xi(n_1)  ,  \xi(n_2) )  \cap    \Span  (  \xi(n'_3)  , \ldots ,  \xi(n'_{p+1}) ) = \{0\}$$
using $(d)$, so that $ \xi(n_1) $,  $\xi(n_2)$, $\xi(n'_3)$, \ldots ,  $\xi(n'_{p+1})$ are linearly independent (using the property analogous to $(b)$ for $n'_3$, \ldots, $n'_{p+1}$). In this case the conclusions of Lemma  \ref{lemalglin} hold  with $n_1$, $n_2$, $n'_3$, \ldots, $n'_{p+1}$ (using $(d)$ and the fact that $n'_3,\ldots, n'_{p+1}\in\calN \setminus\llbracket n_1-\delta, n_2+\delta]$ so that $|n'_i-n_j|>\delta$ for any $i\in\{3,\ldots,p+1\}$ and any $j\in\{1,2\}$).

If $d=1$ then we may assume that  $\xi(n_2)$, $\xi(n''_3)$, \ldots ,  $\xi(n''_{p+2})$ are linearly independent and span a vector space which contains $\xi(n_1)$; indeed otherwise the same properties would hold after permuting $n_1$ and $n_2$. Then the conclusions of Lemma  \ref{lemalglin} hold  with  $n_2$, $n''_3$, \ldots, $n''_{p+2}$.

At last, if $d=2$ then $   \Span  (  \xi(n''_3)  , \ldots ,  \xi(n''_{p+2}) )$ contains both  $\xi(n_1)$ and  $\xi(n_2)$; therefore it contains $\Span (  \xi(n_1) , \ldots,  \xi(n_p) )$ using $(d)$. These vector spaces are therefore equal because they have the same dimension; by construction of $s$, they don't contain $\xi(s)$. Since $s\in\llbracket n_1,n_2\rrbracket $, this is enough to prove that the conclusions of Lemma  \ref{lemalglin} hold  with   $n''_3$, \ldots, $n''_{p+2}$, $s$.

This concludes the proof of Lemma  \ref{lemalglin} in all cases.

\subsection{Proof of Theorem \ref{th70}} \label{secdem70}

 Let $\eps >0$,  and $A\geq D \geq 1$ be such that $0 < \eps \leq 1/20$ and $A\geq \eps^{-12/\eps} D$ (we denote here by capital letters the variables $a$ and $d$ of Theorem \ref{th70}). We choose an odd integer $b$ such that $1+8D/\eps < b < 9D/\eps$, and denote by $a$ the odd integer such that $b$ divides $a$ and $A-3b+2\leq a \leq A-b+1$.  We put $k = (b+1)/2$ and let $\xi(s) = \zeta(2s+1)$ for any $s\in\llbracket 1,(a+b-2)/2\rrbracket $. For any $s \in \llbracket 1,(a-1)/2\rrbracket $ we let also
 
$$v_s = \left( \begin{array}{c} \xi(s) \\   \combitiny{2s+2}2 \xi(s+1) \\   \combitiny{2s+4}4 \xi(s+2 ) \\ \vdots \\  \combitiny{2s+2k-2}{2k-2} \xi(s+k-1)    \end{array} \right) =  \left( \begin{array}{c} \zeta(2s+1) \\   \combitiny{2s+2}2 \zeta(2s+3) \\   \combitiny{2s+4}4 \zeta(2s+5 ) \\ \vdots \\  \combitiny{2s+2k-2}{2k-2} \zeta(2s+b)    \end{array} \right) .$$

Since $A \geq  \eps^{-12/\eps} D  \geq 20^{240} D$, we have $a\geq 9b$ so that Theorem \ref{thdiogen} yields
\begin{equation} \label{eqdioun}
\rk_\Q(e_1,\ldots,e_k,v_1,\ldots,v_{(a-1)/2})\geq  k \sup_{r\in \Iab} \Big( 1- \frac{\log \alrab}{\log \Qrab}\Big)
\end{equation}
where $(e_1,\ldots,e_k)$ is the canonical basis of $\R^k$. Now we let $E = \R/\Q$ and denote by $\pi_0 : \R^k\to E^k$ the canonical surjection on each component. Then Eq. \eqref{eqdioun} yields
\begin{equation} \label{eqdiode}
\rk_\Q( v'_1,\ldots,v'_{(a-1)/2})\geq  k \sup_{r\in \Iab} \Big(  - \frac{\log \alrab}{\log \Qrab}\Big) 
\end{equation}
where $v'_i = \pi_0(v_i)\in E^k$; indeed the restriction of $\pi_0$ to the $\Q$-subspace generated by $e_1,\ldots,e_k,v_1,\ldots,v_{(a-1)/2}$ has kernel equal to $\Q^k$, which has dimension $k$. Letting $\xi'(s)$ denote the image of $\xi(s) = \zeta(2s+1)$ in $E = \R/\Q$ and $\lambda_{i,s} = \combitiny{2s+2i-2}{2i-2} \in \K = \Q $ for $1\leq i \leq k$ and $1 \leq s \leq (a-1)/2$, Proposition \ref{propalglin} applies with $\delta = D/2$, $p=N+1$ and $q=0$ if we use Eq. \eqref{eqdiode} and check that
\begin{equation} \label{eqaverif}
 k \sup_{r\in \Iab} \Big(  - \frac{\log \alrab}{\log \Qrab}\Big)  > (k+2D)N;
 \end{equation}
here $N = [\frac{1-\eps}{1+\log 2} \log(A/D)]$ as in Theorem \ref{th70}.
Assuming (for the time being) that this inequality holds, Proposition \ref{propalglin}  provides integers $n_1, \ldots, n_{N+1}\in\llbracket 1,(a+b-2)/2\rrbracket $ such that
\begin{itemize}
\item $\xi'(n_1)$, \ldots, $\xi'(n_{N+1})$ are $\Q$-linearly independent in $E = \R/\Q$.
\item For any $i,j\in\{1,\ldots,N+1\}$ with $i\neq j$, $|n_i-n_j| > \delta$.
\end{itemize}
We let $\sigma_i = 2n_i+1$ for any $i\in\{1,\ldots,N+1\}$. Then for any $i\neq j$ we have $|\sigma_i-\sigma_j| > 2\delta=D$ so that $\sigma_i > D$ for any $i$ with at most one exception. Reordering the $\sigma_i$'s if necessary, we may assume that $D < \sigma_i \leq a+b-1\leq A$ for any $i\in\{1,\ldots,N\}$. Moreover if 1, $\zeta(\sigma_1)$, \ldots, $\zeta(\sigma_N)$ were linearly dependent over $\Q$, there would exist $\lambda_0,\ldots, \lambda_N\in\Q$, not all zero, such that $\lambda_0 + \lambda_1\zeta(\sigma_1) + \ldots  +  \lambda_N\zeta(\sigma_N)=0$. Seen in the quotient space $E $, this relation reads $ \lambda_1\xi'(n_1) + \ldots  +  \lambda_N\xi'(n_N)=0$. It is non-trivial since $(\lambda_1,\ldots,\lambda_N)\neq(0,\ldots,0)$, so that it contradicts the $\Q$-linear  independence of $\xi'(n_1)$, \ldots, $\xi'(n_{N+1})$  in $E$. Therefore 1, $\zeta(\sigma_1)$, \ldots, $\zeta(\sigma_N)$ are $\Q$-linearly independent real numbers; this concludes the proof of Theorem~\ref{th70}, provided we check Eq. \eqref{eqaverif}.

\bigskip

In order to  check Eq. \eqref{eqaverif}, we recall that $0 < \eps \leq 1/20$ and $A\geq \eps^{-12/\eps} D$, and let $r\geq 1$ denote the integer part of $(A/D)^{1-\eps/3}$. Since the map $x\mapsto x^{-\eps/2}\log(x)$ is non-increasing on $[\exp(2/\eps), +\infty)$ and $A/D \geq\eps^{-12/\eps}  \geq \exp(2/\eps)$, we have
$$\log(A/D)\leq \eps^6 \log( \eps^{-12/\eps} ) (A/D)^{\eps/2} = 12 \eps^5 \log(1/\eps)  (A/D)^{\eps/2}$$
so that
$$br\log(4r+3)< br\log(A/D) \leq 12 b \eps^5 \log(1/\eps)  A/D  < 108 \eps^4  \log(1/\eps)  A < 0.041 \eps A$$
since $b< 9D/\eps$ and $0 < \eps \leq 1/20$. This implies $\frac92 br \log (4r+3)\leq a $ so that $r\in\Iab$, and also
$$4br\log r + 2(a+b-1) + 2b(r+1)\log 2 < 2a+\eps A/4$$
since $b  < 9D/\eps \leq 9 \eps^{ 12/\eps} A /\eps  \leq 9\eps^{239}A$. In the same way we have
$$2(b-1)+2b(2r+1)\log(2r+1)< \eps A / 4.$$
These inequalities yield
\begin{eqnarray*}
 - \frac{\log \alrab}{\log \Qrab}
&=& \frac{ 2(a-2br)\log r  - 2(a+b-1) - 2b(r+1)\log 2}{2(a+b-1)+2(a-2br)\log 2 + 2b(2r+1)\log(2r+1)}\\
&>& \frac{2a\log r -2a - \eps A/4}{2(1+\log 2)a + \eps A/4} > \frac{\log r -1-\eps/7.99}{1+\log 2 + \eps/7.99}
\end{eqnarray*}
since $A \leq a+3b\leq a + 27 \eps^{239}A \leq 8a/7.99.$ Moreover we have $r \geq (A/D)^{1/2}\geq \eps^{-6/\eps}$ so that
$$ \log r -1-\eps/7.99 \geq \log(r+1) - \frac1r-1-\eps/7.99 \geq (1-\eps/3)\log(A/D)-1-\eps/7.9.$$
On the other hand,  
\begin{eqnarray*}
\frac{k+2D}{k}N 
&=& \Big(1+4D/(b+1)\Big)\Big[\frac{1-\eps}{1+\log 2}\log(A/D)\Big]\\
&\leq&1+\eps/2 + \frac{(1+\eps/2)(1-\eps)}{1+\log 2} \log(A/D).
\end{eqnarray*}
Combining these inequalities yields
$$
 - \frac{\log \alrab}{\log \Qrab} - \frac{k+2D}{k}N  > \frac{g(\eps)}{1+\log 2} \log(A/D) - h(\eps)$$
 with 
 $$g(\eps) = \frac{1-\eps/3}{1+\frac{\eps}{7.99(1+\log 2)}} - (1+\eps/2)(1-\eps) \geq 0.09 \eps$$
 and 
 $$h(\eps) = \frac{1+ \eps/7.9}{1+\log 2} + 1 +\eps/2 \leq 1.62$$
 since $\eps \leq 1/20$. Now $\eps \log(A/D) \geq 12\log 20 \geq 35$ so that
 $$
 - \frac{\log \alrab}{\log \Qrab} - \frac{k+2D}{k}N  > \frac{0.09}{1+\log 2}\cdot 35 - 1.62 \geq 0.24 > 0.$$
 This concludes the proof of Eq. \eqref{eqaverif}, and that of Theorem \ref{th70}.

\subsection{Proof of Theorem \ref{thzetafacile}}  \label{subsecdemzetafacile}

To prove Theorem \ref{thzetafacile}, we follow the proof of Theorem \ref{th70} and let $\pi : \R^{(b+1)/2}\to\R$ be defined by $\pi(x_1,x_2,\ldots , x_k) = \lambda_0 x_1+\lambda_1x_2+\ldots+\lambda_d x_{d+1}$. We apply Remark \ref{Rk2} that follows Theorem \ref{thdiogen}, so that Eq. \eqref{eqdioun} is replaced with 
$$
\rk_\Q(1,\pi(v_1),\ldots,\pi(v_{(a-1)/2}))\geq  \sup_{r\in \Iab} \Big( 1- \frac{\log \alrab}{\log \Qrab}\Big).
$$
Since $\pi(v_{(s-1)/2})$ is exactly the number \eqref{eqfac}, Theorem \ref{thzetafacile} follows from the fact that this lower bound is greater than $N+1$ (see Eq. \eqref{eqaverif}).

\subsection{Proof of Theorem \ref{th145}} \label{secdem145}

In this section we prove Theorem \ref{th145} stated in the introduction, by following the proof of Theorems \ref{thdiogen} and \ref{th70} (see \S\S \ref{seczeta} and  \ref{secdem70}). We indicate simply the differences.

Denoting by $D$ the odd integer $d$ in the statement of Theorem \ref{th145}, we let $b = D$ if $ D \leq 20000$, and otherwise we define $b$  to be  the least odd integer $b$ such that $D \leq \csi b$. We put also $a = 149 b$, $r=11$, and $k = (b+1)/2$. We have $\frac92 br\log(4r+3)\leq a$, and we may assume that $r\not\in\Rab$ (otherwise we replace everywhere $r$ with $r+\eps$ for a sufficiently small rational number $\eps>0$). The real root $\mu_1$ of \S \ref{subsecasydebut} is independent from $b$, since it is the unique root of the polynomial 
$ (X+23)  (X-1)^{150} -  (X-23)  (X+1)^{150}$ in the real interval $(23,+\infty)$; we have $\mu_1 \simeq \dun$. Since $f'(\mu_1+i0) = bi\pi\in i\R$, we obtain using Eq.~\eqref{eq212}:
$$\eps_{ b} =   \exp \Re  f_0(\mu_1+i0)=   \exp \Re  f (\mu_1+i0)  \simeq \exp(-\dde b).$$
We shall use this numerical value instead of the last upper bound of Eq. \eqref{eq316bis}. 

Another refinement turns out to be necessary to complete the proof with these parameters: Zudilin has constructed 
in Proposition 1 of \cite{Zudilincentqc} a sequence $(\Pi_n)_{n\geq 1}$ of positive integers such that $\Pi_n^{-b}L_n$ is still a linear form with integer coefficients, and 
$$\lim_{n\to\infty} \Pi_n^{1/n} = \varpi \simeq \dtr.$$

If $ D \leq 20000$, we let $t=1$ and notice that the lower bound on $\sigma_2-\sigma_1$ in Theorem \ref{th145} is equivalent to $\sigma_2 > \sigma_1$. Otherwise, we let $t$ denote the integer part of $ \cun b $, and remark (for future reference) that $150b \leq  151 D$.  We  define $\pi : \R^{(b+1)/2}\to\R^t$ by $\pi(x_1,x_3,\ldots,x_b) = (x_{b-2t+2},\ldots,x_{b-2},x_b)$. Following the proof of Theorem \ref{thdiogen} and Remark~\ref{Rk2}, we obtain
$$\rk_\Q(u'_1,u'_2,\ldots,u'_t,\pi(v_1),\pi(v_2) , \ldots,  \pi(v_{(a-1)/2})) \geq t  \Big( 1 - \frac{\log \alpha}{\log Q}\Big) $$
where $(u'_1,\ldots,u'_t)$ is the canonical basis of $\R^t$, $v_1,\ldots,v_{(a-1)/2}$ are defined in Theorem \ref{thdiogen}, and 
$$\alpha = \exp(2(a+b-1)-b\varpi+\Re  f (\mu_1+i0))\simeq \exp(- \cde b-2),$$
$$Q = \exp(2(a+b-1)-b\varpi+2(a-2br)\log 2+2b(2r+1)\log(2r+1))\simeq \exp(\ctr b -2).$$
Now we let $E = \R/\Q$  and denote by $\pi_0 : \R^t\to E^t$ the canonical surjection on each component, as in \S \ref{secdem70}. 
 Letting also $v'_i = \pi_0(\pi(v_i))$ we obtain
$$\rk_\Q( v'_1  , \ldots,   v'_{(a-1)/2} ) \geq  - t  \, \frac{\log \alpha}{\log Q} > t + 4\delta$$
 where  $  \delta = \cqu b $  if  $ D > 20000$,   and $\delta = 0$  otherwise.  
Let us  consider $\xi(n ) = \zeta(b-2t+2n+2)$ for $1\leq n \leq N$, with $N = \frac{a-3}2+t$. Proposition \ref{propalglin}, applied with $p=2$ and $q=0$, provides integers $n_1,n_2\leq N$ such that $n_2 > n_1+\delta$. Letting $\sigma_i = b-2t+2n_i+2$ for $i\in\{1,2\}$, we obtain if $D > 20000$:
$$\sigma_2 > \sigma_1 + 2\delta =  \sigma_1 + \deuxcqu b \geq \sigma_1+\cci D,$$
$$D+2 \leq \csi b+2 \leq b-2t+4\leq\sigma_1<\sigma_2\leq b+a-1 = 150b-1\leq 151 D,$$
and 1, $\zeta(\sigma_1)$, $\zeta(\sigma_2)$ are $\Q$-linearly independent. If $D\leq 20000$ the last inequality is simply replaced with
$$D+2 = b+2= b-2t+4\leq\sigma_1<\sigma_2\leq b+a-1 = 150b-1\leq 150 D.$$
In both cases this concludes the proof of Theorem \ref{th145}.

\bigskip

Let us conclude this section with two remarks on the proof. First, taking $b=d$ and $t=1$, the proof yields 
$$\dim_\Q\Span_\Q(1,\zeta(d+2),\zeta(d+4),\ldots,\zeta(150d-1)) \geq 3$$
for any $d\geq 1$. On the other hand, if $b$ is small then the estimates in the proof are slightly better (see the definition of $\alpha $ and $Q$); this is the reason why Zudilin obtains 145 instead of 151 when $b=1$. The improvements of \cite{SFZu}, leading to the value 139, fall also into error terms as $b\to\infty$.

\subsection{Proof of Corollary \ref{cor82}} \label{secdemcor82}

In this section we deduce Corollary \ref{cor82} from Theorem \ref{th70}.

\bigskip

Let $\eps > 0$ be such that $\eps \leq 1/20$; put $\eta = \eps^{15/\eps}$ and
$$\eps' = \eps - \frac{\eps}{3\log(1/\eps)} $$
so that $0 < \eps' < \eps \leq 1/20$. We consider also
$$\Chd = {\eps'}^{-12/\eps'}, M = \Big[ \frac{1-\eps'}{1+\log 2}\log \Chd\Big], \mbox{ and } \Ctihd = (2e)^{(1+\eps)M}.$$
We shall prove below that
\begin{equation} \label{eqaverifde}
 \Chd+1 \leq \Ctihd <  \eta^{-1}.
 \end{equation}
 Taking this inequality for granted, we assume (by induction on $n$) that $u_1$, $u_2$, \ldots, $u_{Mn}$ are already defined so that the first three conclusions of Corollary \ref{cor82} hold for any $i \leq Mn$, with $u_{Mn}\leq \Ctihd^{n-1} \Chd$. If $n=0$ this assumption is empty. Then we apply Theorem \ref{th70} with $\eps'$, $a =  \Ctihd^n \Chd$ and $d =  \Ctihd^n $ (since $a/d = \Chd =  {\eps'}^{-12/\eps'}$ and $0 <  \eps'  \leq 1/20$). This provides odd integers $u_{Mn+1}$, \ldots, $u_{M(n+1)}$ such that $ \Ctihd^n < u_{Mn+1} <  \ldots < u_{M(n+1)} \leq  \Ctihd^n \Chd$, $\zeta(u_{Mn+j})\not\in\Q$ for any $j\in\unM$, and $u_{Mn+j+1}-u_{Mn+j}> \Ctihd^n$ for any $j\in\unMmu$. Using Eq. \eqref{eqaverifde} this lower bound implies
 $$u_{Mn+j+1} / u_{Mn+j} > 1 +  \Ctihd^n / u_{Mn+j} \geq 1+1/ \Chd \geq 1+\eta.$$
For $j=0$ we obtain in the same way the following inequalities, since $u_{Mn} \leq \Ctihd^{n-1} \Chd$:
 $$u_{Mn +1} / u_{Mn } > \frac{ \Ctihd^n }{ \Ctihd^{n-1}\Chd} \geq  \Ctihd / \Chd \geq (\Chd+1)/\Chd \geq 1+\eta.$$
Letting $i=Mn+j$ with $1 \leq j  \leq M$, we have $i/M-1\leq n < i/M$ so that, using  Eq. \eqref{eqaverifde} again:
$$\eta (2e)^{(1+\eps)i}  \leq \Ctihd^{-1} \eta (2e)^{(1+\eps)i} =  \Ctihd^{i/M-1} \leq \Ctihd^n < s_i \leq  \Ctihd^n \Chd < \Ctihd^{i/M} \Chd  \leq
 \eta^{-1} (2e)^{(1+\eps)i} .$$
 This concludes the induction.
At last, if  $a \geq \eta^{-1/\eps}$ then letting $N = [\frac{1-2\eps}{1+\log 2}\log a]$ this upper bound on $s_N$ yields
$$s_N < \eta^{-1} a ^{1-\eps-2\eps^2}  < a\eta^{2\eps} < a.$$
This concludes the proof of Corollary \ref{cor82}, except for Eq.  \eqref{eqaverifde} that we shall prove now.

\bigskip

To begin with, we notice that $3\log(1/\eps) \geq 3\log (20) \geq 8.98$ so that $\eps' \geq 0.88 \eps$ and $\eps \leq 1.14 \eps'$. This implies
$\log(1/\eps) \geq \log(1/\eps') + \log(0.88)> 19 \log(1/\eps')/20$ since $\log(20)/20 > -  \log(0.88)$, and finally:
\begin{equation} \label{eqepsquo}
\eps' / \log(1/\eps') < \eps / \log(1/\eps)  < 1.2 \eps' / \log(1/\eps').
  \end{equation}
This enables us to prove that $\Ctihd <  \eta^{-1}$ because
\begin{eqnarray*}
\log \Ctihd 
&\leq& (1+\eps)(1-\eps')\log \Chd \leq \Big(1+\frac{\eps}{3\log(1/\eps)}\Big)\log \Chd\\
&\leq& 1.006 \log \Chd \mbox{ since } \eps\leq 1/20 \mbox{ and } 3\log(1/\eps) \geq 3 \log 20 \geq 8.98\\
&<&  \log(\eta^{-1}) \mbox{ by definition of $\Chd$ and $\eta$, and using Eq. \eqref{eqepsquo}.}
\end{eqnarray*}
At last, we have
\begin{eqnarray*}
\log \Ctihd  - \log(\Chd+1)
&\geq& (1+\eps)N(1+\log2) - \log \Chd  - \log(1+1/\Chd)\\
&\geq& [(1+\eps) (1-\eps')-1]\log \Chd -  (1+\eps) (1+\log2) - 1/\Chd\\
&\geq& 12 (\eps/\eps' -1- \eps)\log(1/\eps') -  (1+\eps) (1+\log2) -\eps'^{12/\eps'}.
\end{eqnarray*}
Since $\eps/\eps' = (1-1/(3\log(1/\eps)))^{-1} \geq 1+1/(3\log(1/\eps))$ and $\log(1/\eps') \geq \log(1/\eps )$ we obtain
$$\log \Ctihd  - \log(\Chd+1) \geq 4 - 12\eps \log(1/\eps ) - (1+\eps)(1+\log 2) -\eps'^{12/\eps'} > 0.4 > 0$$
since $\eps' < \eps  \leq 1/20$.
This concludes  the proof of  Eq.  \eqref{eqaverifde},  and that of Corollary \ref{cor82}.

\newcommand{\url}{\texttt}

\providecommand{\bysame}{\leavevmode ---\ }
\providecommand{\og}{``}
\providecommand{\fg}{''}
\providecommand{\smfandname}{\&}
\providecommand{\smfedsname}{\'eds.}
\providecommand{\smfedname}{\'ed.}
\providecommand{\smfmastersthesisname}{M\'emoire}
\providecommand{\smfphdthesisname}{Th\`ese}

\end{document}